\documentclass[10pt,a4paper,twoside]{article}

\usepackage{a4}
\usepackage{exscale}
\usepackage[]{amsmath}
\usepackage{amssymb}
\usepackage{amsfonts}
\usepackage{eucal}
\usepackage[amsmath,thmmarks]{ntheorem}
\usepackage{bbm}

\theoremseparator{:}
\newtheorem{theorem}{Theorem}[section]
\newtheorem{lemma}[theorem]{Lemma}

\newtheorem{proposition}[theorem]{Proposition}
\newtheorem{corollary}[theorem]{Corollary}

\theorembodyfont{\normalfont}
\newtheorem{remark}[theorem]{Remark}
\newtheorem*{remarkon}{Remark}
\newtheorem{remarks}[theorem]{Remarks}

\theoremsymbol{\ensuremath{\square}}
\theoremheaderfont{\normalfont\scshape}
\theorembodyfont{\normalfont}
\newtheorem*{proof}{Proof}

\newcommand{\R}{\ensuremath{\mathbb{R}}}
\newcommand{\N}{\ensuremath{\mathbb{N}}}

\newcommand{\ccdot}{\,\cdot\,}

\renewcommand{\epsilon}{\varepsilon}

\usepackage[hang,small,bf]{caption}
\usepackage{fancyhdr}
\usepackage[active]{srcltx}

\usepackage[top=3.0cm,bottom=2.6cm,left=2.6cm,right=2.6cm]{geometry}
\theoremlisttype{all}

\DeclareMathOperator{\Sing}{Sing}
\DeclareMathOperator{\Reg}{Reg}

\DeclareMathOperator{\diverg}{div}
\DeclareMathOperator{\dist}{dist}
\DeclareMathOperator{\diam}{diam}
\DeclareMathOperator{\supp}{spt}

\newcommand{\Hm}{\mathcal{H}}

\newcommand{\Lm}{\mathcal{L}}

\newcommand{\Lk}{\mathcal{L}^k}
\newcommand{\dx}{\, dx}
\newcommand{\dy}{\, dy}
\newcommand{\dt}{\, dt}

\newcommand{\kabs}[1]{\ensuremath{\vert#1\vert}}
\newcommand{\babs}[1]{\ensuremath{\big\vert#1\big\vert}}
\newcommand{\Babs}[1]{\ensuremath{\Big\vert#1\Big\vert}}

\newcommand{\knorm}[1]{\ensuremath{\Vert#1\Vert}}

\newcommand{\tmint}{\mathop{\int\hskip -0,89em -\,}}
\newcommand{\mint} {\mathop{\int\hskip -1,05em -\,}}

\newcommand{\mI}[1]{\mint\nolimits_{\!\!\!\!#1}}

\usepackage{titlesec}  
\titleformat{\section}{\bf\large}{\thesection}{1em}{}
\titleformat{\subsection}{\bf\large}{\thesubsection}{1em}{}
\titleformat{\subsubsection}{\bf\normalsize}{\thesubsubsection}{1em}{}

\numberwithin{equation}{section}
\allowdisplaybreaks[3] 

\title{
{\bf Boundary regularity for elliptic systems under a natural growth condition}
}
\date{}
\author{Lisa Beck\footnote{L. Beck, Scuola Normale Superiore di Pisa, Piazza dei Cavalieri 7, Pisa, Italy. E-mail: lisa.beck@sns.it}}

\begin{document}
 

\maketitle

\begin{abstract}
We consider weak solutions $u \in u_0 + W^{1,2}_0(\Omega,\R^N) \cap L^{\infty}(\Omega,\R^N)$
of second order nonlinear elliptic systems of the type
\begin{equation*}
- \diverg  a (\, \cdot \, , u, Du ) \, = \, b(\, \cdot \, ,u,Du) \qquad \text{ in } \Omega
\end{equation*}
with an inhomogeneity satisfying a natural growth condition. In dimensions $n \in \{2,3,4\}$ we show that $\Hm^{n-1}$-almost every boundary point is a regular point for $Du$, provided that the boundary 
data and the coefficients are sufficiently smooth.  \vspace{0.2cm}

\emph{Mathematics Subject Classification (2000):} 35J45, 35J55

\end{abstract}

\section{Introduction}

In this paper we are concerned with the existence of regular boundary points for the gradient of bounded, vector-valued weak solutions $u \in W^{1,2}(\Omega,\R^N) \cap L^{\infty}(\Omega,\R^N)$ of nonlinear, inhomogeneous elliptic systems of the form
\begin{equation}
\label{DP-rbp}
-\diverg a (\, \cdot \, , u, Du ) \, = \,  b(\, \cdot \, ,u,Du)  \qquad  \text{in } \Omega
\end{equation}
with boundary values $u_0$ on $\partial \Omega$ in the sense of traces. Here $\Omega \subset \R^n$ is a bounded domain of class $C^{1,\alpha}$ and $u_0 \in C^{1,\alpha}(\overline{\Omega},\R^N)$ for some $\alpha \in (0,1)$. The coefficients $a:\Omega \times \R^N \times \R^{nN} \to \R^{nN}$ are assumed to be H\"older continuous with exponent $\alpha$ with respect to the first two variables and of class $C^1$ in the last variable, satisfying a standard quadratic growth condition. Furthermore, we assume that the right-hand side $b \colon \Omega \times \R^N \times \R^{nN} \to \R^N$ satisfies a natural growth condition and that an additional smallness condition on $\knorm{u}_{L^\infty}$ holds. In general we cannot expect a weak solution to a nonlinear elliptic system -- in contrast to weak solutions to a single equation~-- to be a classical one of class $C^2$, see  \cite{DEGIORGI68,GIUSMIRA68EX}. Nevertheless, a partial regularity result still holds true which can be stated as follows: every weak solution $u \in W^{1,2}(\Omega,\R^N) \cap L^{\infty}(\Omega,\R^N)$ to the inhomogeneous system \eqref{DP-rbp} is of class $C^1$ near a point $x_0$ if and only if a certain excess quantity is sufficiently small and the mean values of $u$ and of $Du$ on balls $B_{\rho}(x_0)$ do not not diverge for $\rho \searrow 0$. To be more precise, if we denote by
\begin{equation*}
\Reg_{Du}(\overline{\Omega}) \, := \, \big\{ x_0 \in \overline{\Omega}: Du \in C^0(\overline{\Omega} \cap A,\R^{nN})  \text{ for some neighborhood } A \text{ of } x_0 \big\} 
\end{equation*}
the set of regular points for $Du$ (in the interior and at the boundary), and by $\Sing_{Du}(\overline{\Omega}) := \overline{\Omega} \setminus \Reg_{Du}(\overline{\Omega})$ the set of singular points of $Du$, then the singular set is characterized via $\Sing_{Du}(\overline{\Omega}) = \Sigma \cup \Sigma_u$ with
\begin{gather*}
\Sigma \, := \, \Big\{ x_0 \in \overline{\Omega} \colon \liminf_{\rho \to \, 0^+} \mI{\Omega \cap B_{\rho}(x_0)} \!\!\! \big\vert Du - (Du)_{\Omega \cap B_{\rho}(x_0)}\big\vert^2 \dx > 0 
	\text{ or } \limsup_{\rho \to \, 0^+} 
	\babs{(Du)_{\Omega \cap B_{\rho}(x_0)}} = \infty \Big\} \,, \\
\Sigma_u \, :=  \, \Big\{ x_0 \in \overline{\Omega} \colon \limsup_{\rho \to \, 0^+} 
	\babs{(u)_{\Omega \cap B_{\rho}(x_0)}} = \infty \Big \}\,.
\end{gather*}
Furthermore, the gradient $Du$ of the weak solution is locally H\"older continuous with exponent $\alpha$ in a (small) neighborhood of every point $x_0 \in \Reg_{Du}(\overline{\Omega})$, see \cite{GROTOWSKI02} (and \cite{HAMBURGER07,BECK07} for the non-quadratic analogues). This is the up-to-the-boundary extension of the interior partial regularity results obtained in various papers starting from \cite{GIUSMIRA68,GIAMOD79,IVERT79}. By Lebesgue's differentiation theorem, the regularity criterion stated above applies to almost every point in $\overline{\Omega}$, whence $\kabs{ \overline{\Omega} \setminus \Reg_{Du}(\overline{\Omega})} = 0$. However, this does not yield the existence of even one single regular boundary point for weak solutions to general nonlinear elliptic systems, since the boundary $\partial \Omega$ itself is a set of Lebesgue measure zero. By contrast, due to Giaquinta's counterexample \cite{GIAQUINTA78}, it is well known that singularities may occur at the boundary even if the boundary data is smooth.

The question of dimension reduction of the singular set (in the sense that it is not only negligible with respect to the Lebesgue measure but that its Hausdorff dimension is bounded strictly below $n$) has received considerable attention in recent years. Some significant results were first obtained for weak solutions to systems satisfying special structure conditions: for quasilinear systems of the form
\begin{equation*}
- \diverg \big( a(\ccdot,u) \, Du \big) \, = \, b(\ccdot,u,Du) \,,
\end{equation*}
various partial regularity results were established, stating that the weak solution $u$ (instead of its first derivative) is locally H\"older continuous. To bound the Hausdorff dimension of the singular set $\Sing_{u}(\overline{\Omega})$, we recall that a regular point $x_0 \in \overline{\Omega}$ of $u$ is a point where $u$ is locally continuous and is characterized via a smallness condition on the lower order excess functional
\begin{equation*} 
\mI{\Omega \cap B_{\rho}(x_0)} \big\vert u - (u)_{\Omega \cap B_{\rho}(x_0)}\big\vert^2 \dx \,,
\end{equation*}
e.\,g., see \cite{GIUSMIRA68,COLOMBINI71,PEPE71,GROTOWSKI02QL,ARKHIPOVA96}. Since the set of non-Lebesgue points of every $W^{1,p}$-map has Hausdorff dimension not larger than $n-p$, the Hausdorff dimension of $\Sing_{u}(\overline{\Omega})$ cannot exceed $n-2$. If the coefficient matrix $a(\cdot,\cdot)$ of the quasilinear system is further assumed to be of diagonal form, it is known that the weak solution is a classical solution (see \cite{WIEGNER76} where boundary regularity is included). Useful estimates for the singular set are also available for nonlinear elliptic systems obeying special structure assumptions: for instance, Uhlenbeck established in her fundamental paper \cite{UHLENBECK77} a strong maximum principle for the gradient $Du$ of weak solutions to nonlinear systems, provided that the nonlinear part of the coefficient function only depends on the modulus of $Du$. This was the key to obtain everywhere-regularity for $Du$. For an extension to the nonquadratic case we refer to \cite{TOLK83,ACEFUS89}. However, neither could Uhlenbeck's techniques be carried over to the boundary, nor is a suitable counterexample available in the literature, leaving the question of full boundary regularity open for such systems. Turning the attention to general nonlinear elliptic systems, we observe that a direct comparison technique allows to infer local H\"older continuity of the weak solution outside a set of Hausdorff dimension $n-p$ in low dimensions $n \leq p+2$, see  \cite{CAMPANATO82,CAMPANATO87,ARKHIPOVA97,ARKHIPOVA03,BECK09b}. By contrast, in \emph{arbitrary} dimensions $n$ the reduction of the Hausdorff dimension of the singular set $\Sing_{Du}(\overline{\Omega})$ for the \emph{gradient} $Du$ was a long-standing problem. It was finally tackled by Mingione \cite{MINGIONE03arch}: he studied the interior singular set $\Sing_{Du}(\Omega)$ in the superquadratic case $p \geq 2$ for systems without $u$-dependencies and with inhomogeneities obeying a controllable growth condition, and he succeeded in showing that the Hausdorff dimension of $\Sing_{Du}(\Omega)$ is not larger than $n - 2 \alpha$. In \cite{MINGIONE03} he extended these results to systems with inhomogeneities under a natural growth condition, covering also systems explicitly depending on $u$, provided that $n \leq p+2$ is satisfied. 

We now return to the existence of regular boundary points: we first observe that for this aim the almost-everywhere regularity result has to be improved to a bound for the Hausdorff dimension less than $n-1$ because this yields immediately that almost every boundary point is regular. Consequently, our objective is to identify additional assumptions on the coefficients or on the space dimension which guarantee this dimension reduction. A result in this direction was recently obtained by Duzaar, Kristensen and Mingione \cite{DUZKRIMIN05}: they considered weak solutions $u \in W^{1,p}(\Omega,\R^N)$, $p \in (1,\infty)$, of the homogeneous Dirichlet problem corresponding to \eqref{DP-rbp} and developed a technique which allows to carry the estimates in \cite{MINGIONE03arch} up to the boundary, implying in particular the existence of regular boundary points, provided that $n-2\alpha < n-1$ (or equivalently $\alpha > \frac{1}{2}$) is satisfied. More precisely, the authors obtained for every $\alpha \in (\frac{1}{2},1]$ that almost every boundary point is regular if the coefficients $a(x,z)$ have no $u$-dependency or if $n \leq p+2$ holds. In the quadratic case they improved this result in two ways: on the one hand, inhomogeneities with controllable growth were included, and on the other hand the condition on $\alpha$ was sharpened to $\alpha > \frac{1}{2} - \epsilon$ for some number $\epsilon > 0$ stemming from an application of Gehring's lemma. We further mention that various results establishing better estimates for the (interior) singular set of minimizers of variational integral can be found in \cite{KRIMIN05,KRIMIN08}. 

The main result in this paper is an extension of the result \cite{DUZKRIMIN05} to bounded weak solutions to inhomogeneous systems under a critical growth condition on the inhomogeneity (giving also an alternative proof of \cite[Theorem 1.3]{DUZKRIMIN05}), namely the improvement of the estimate $\kabs{\Sing_{Du}(\overline{\Omega})}=0$ in the following sense:

\begin{theorem}
\label{regul-bp-mit-u-nat-2}
Consider $n \in \{2,3,4\}$ and $\alpha \geq \frac{1}{2}$. Let $\Omega \subset \R^n$ be a domain of class $C^{1,\alpha}$ and $u_0 \in C^{1,\alpha}(\overline{\Omega},\R^N)$. Assume that $u \in u_0 + W^{1,2}_0(\Omega,\R^N) \cap L^{\infty}(\Omega,\R^N)$ is a weak solution of the Dirichlet problem \eqref{DP-rbp} under the assumptions (H1)-(H3) and (B) from Section~\ref{struc-not}, and suppose that $\knorm{u}_{L^{\infty}(\Omega,\R^N)} \leq M$ for some $M>0$ such that $2 L_2 M < \nu$. Then $\Hm^{n-1}$-almost every boundary point is a regular point for $Du$.
\end{theorem}

This result was presented as a part of the author's PhD thesis \cite{BECK08} where most of the proofs and calculations are discussed in detail. In addition, some extensions and open questions concerning the dimension reduction of the singular set are collected in Section~\ref{sec-ext}. In particular, in view of an observation by Kristensen and Mingione \cite{KRIMIN08}, it is possible to replace a part of the H\"older continuity assumption on the coefficients with respect to the $x$-variable (in the sense that it is only required with an arbitrary exponent) by an additional fractional differentiability assumption on the map $x \mapsto a(x,u,z)$, see Theorem~\ref{thm-frac}.

We close this introductory part with some remarks about the ideas behind the arguments and the techniques used within this paper. The strategy can be described as follows: To simplify matters we initially consider coefficients of the form $a(x,z)$: If they are H\"older continuous in $x$ with arbitrarily small exponent, we know $\dim_{\Hm}(\Sing_{Du}(\overline{\Omega})) \leq n$. If they are instead Lipschitz-continuous, then standard difference quotients reveal $Du \in W^{1,2}(\Omega,\R^{nN})$ which implies that $\dim_{\Hm}(\Sing_{Du}(\overline{\Omega})) \leq n-2$. Therefore, the upper bound on the Hausdorff dimension of $\Sing_{Du}(\overline{\Omega})$ reflects the regularity of the coefficients in $x$. This gives the impression that the regularity of the coefficients is related not only to the regularity of the solution (namely the local H\"older continuity of $Du$ to the same exponent), but also to the size of the singular set. Working from this observation, Mingione \cite{MINGIONE03arch,MINGIONE03} introduced a remarkable new technique and accomplished in the interior an interpolation between Lipschitz continuity on the one hand and H\"older continuity on the other: for general $\alpha$-H\"older continuous coefficients the existence of higher order derivatives of the weak solution cannot be expected, but it is still possible to differentiate the system \eqref{DP-rbp} in a fractional sense. This leads to the desired upper bound $n - 2 \alpha$ for the Hausdorff dimension. If the coefficients $a(x,u,z)$ now depend explicitly on $u$, the situation becomes more complex and the estimates are technically much more involved. To follow the line of arguments above we have to investigate the regularity of the map $x \mapsto (x,u(x))$. If the weak solution $u$ is a~priori known to be everywhere H\"older continuous then $x \mapsto (x,u(x))$ is also H\"older continuous and the arguments apply with only marginal modifications. However, in general this map is no longer continuous, because $u$ may exhibit irregularities. Nevertheless, at least in low dimensions $n \leq p+2$, local H\"older continuity of weak solutions is guaranteed outside of closed subsets of Hausdorff dimension less than $n-p$. In other words, the set of points where $u$ is not continuous -- and where  $x \mapsto (x,u(x))$ is not regular~-- has sufficiently small Hausdorff dimension, hence, restricting the analysis of $Du$ to the regular set $\Reg_{u}(\Omega)$ of $u$, we still arrive at a good result for $\dim_{\Hm}(\Sing_{Du}(\Omega))$, see Theorem~\ref{int-Mingione}.

In the interior this method relies essentially upon finite difference operators, fractional differentiability estimates for the gradient $Du$ and interpolation techniques dating back to Campanato \cite{CAMPCAN81,CAMPANATO82a}, combined in a delicate iteration scheme (applied for elliptic and parabolic systems \cite{MINGIONE03arch,MINGIONE03,DUMIN05}), and the necessary estimates are deduced by testing with (differences of) the solution. At the boundary some severe problems are caused by the fact that testing is allowed only for differences in tangential direction: hence, the normal direction still has to be recovered by exploiting the system of equations (which follows immediately if second-order derivatives exist). This problem was overcome first for homogeneous elliptic system (and inhomogeneous systems under controllable growth) by an indirect approach developed by Duzaar, Kristensen and Mingione \cite{DUZKRIMIN05}: via a regularization procedure involving both the original coefficients $a(\cdot,\cdot,\cdot)$ and the specific solution $u$, a family of comparison maps is constructed for which the existence of second-order derivatives is known. This allows to gain higher integrability for $Du$ which in turn is used to improve the integrability of the comparison map by means of Calder\'on-Zygmund estimates (provided in \cite{KRIMIN05,BECK08}) in the next iterative step. 

When trying to apply this approach for inhomogeneous systems under critical growth, several critical difficulties arise: most importantly the propagation of higher integrability via the Calder\'on-Zygmund theory seems not to be clear since the natural growth condition merely gives $L^{1+\delta}$ for the right-hand side with some (small) $\delta > 0$ (coming from the higher integrability of $Du$) rather than the necessary prerequisite $L^{q/(p-1)}$ for some $q>p$. For this reason we exploit the system differently and replace the indirect comparison principle by a direct method, introduced by Kronz \cite{KRONZHABIL} as a promising approach for up-to-the-boundary regularity results including upper bounds for the Hausdorff dimension of the singular set, with the flexibility to attack higher order systems. Kronz observed that estimates for the tangential differences suffice to control the averaged mean deviation with respect to mean values taken over slices in tangential direction. Using an alternative definition of fractional Sobolev spaces based on pointwise inequalities, this helps to deduce a fractional differentiability property for the system, which then gives further information on the gradient of the solution. The overall strategy remains unchanged, i.\,e. existence of regular boundary points is still proved by a dimension reduction argument for the singular set $\Sing_{Du}(\overline{\Omega})$: The key tool here is the observation that if $Du$ belongs to a fractional Sobolev space $W^{\theta,p}$, then the characterization of $\Sing_{Du}(\overline{\Omega})$ and a measure density result allow to conclude that the Hausdorff dimension of $\Sing_{Du}(\overline{\Omega})$ does not exceed $n-\theta p$. The proof of such a fractional differentiability estimate for $Du$ is now sketched in a series of steps:

\subsection*{Strategy of the proof:}

\emph{Simplifications:} It suffices to consider the model situation $\Omega = B^+$ and solutions $u \in W^{1,2}(B^+,\R^N) \cap L^{\infty}(B^+,\R^N)$ which vanish on the flat part of the boundary. The general situation then follows from a transformation argument. Furthermore, we assume H\"older continuity of $u$ on $B^+$ with some exponent $\lambda>0$. This is justified by the fact that the solution is H\"older continuous outside a set of Hausdorff dimension $n-2$ in dimensions $n \in \{2,3,4\}$.

\emph{Tangential differences:} Testing the system with differences of the solution up to the boundary is only allowed for tangential directions (because zero boundary values on the flat part are maintained for the test function). Taking into account the assumptions on the coefficients and the inhomogeneity, we end up with an integral estimate for $\kabs{Du(x+h e_s) - Du(x)}$ for all unit directions $e_s \perp e_n$, telling that its $L^2$-norm decays like $c \kabs{h}^{\alpha \lambda / 2}$ (with $\alpha$ denoting the H\"older exponent of the continuity condition on the coefficients with respect to the first and the second variable).

\emph{An estimate for tangential derivatives:} If finite differences of the full derivative $Du$ are estimated, then it is reasonable that also normal differences of only the tangential derivative denoted by $D'u$ are estimated similarly (if we think of Lipschitz-continuous coefficients $a(x,Du)$ for example, this observation is trivial since the previous step yields the existence of second order derivatives $D'Du = D D'u$). This is in fact true (up to a small loss in the power of $\kabs{h}$), and we thus get a first fractional differentiability estimate for $D'u$.

\emph{Towards the normal derivative:} Information about $D_nu$ can only be gained out of the system (in case of Lipschitz-continuous coefficients $a(x,Du)$, the existence of the second order normal derivative $D_{n} Du$ is obtained from the system of equations $- D_n a_n(x,Du) = \sum_{k=1}^{n-1} D_k a_k(x,Du) + b(x,u,Du)$ in a standard way). Looking at the simple example $- \diverg \big( f(x,u) Du \big) = b(x,u,Du)$
we get a first idea on how the coefficients might serve to improve the differentiability of $D_nu$, because we then have $a_n(x,u,Du) = f(x,u) D_nu$, meaning that $a_n(x,u,Du)$ is the missing normal derivative up to a H\"older continuous perturbation (a similar property holds true for the general coefficients). For the moment let us concentrate on $a_n(x,u,Du)$: mimicking the differentiable situation to a certain extent, we show by means of the estimates for tangential differences of $Du$ that slice-wise mean values of $a_n(x,u(x),Du(x))$ are differentiable in the weak sense in the $x_n$-direction, and as a consequence, we obtain that the map $x \mapsto a_n(x,u(x),Du(x))$ is in a fractional Sobolev space.

\emph{An estimate for the normal derivative:} Taking advantage of the ellipticity and the boundedness condition assumed for the coefficients, we find that differences of $D_n u$ are essentially dominated by those of $a_n(x,u,Du)$ and of the tangential derivative $D'u$. Together with a corresponding estimate for the tangential derivatives of $u$, this leads to a fractional differentiability result for the full gradient $Du \in W^{\alpha \gamma \lambda,2}$ for every $\gamma < 1$. 

\emph{Getting rid of $\lambda$:} By an interpolation technique we gain higher integrability out of the fractional differentiability of $Du$. This in turn is used to improve the differentiability of $Du$ in a suitable iteration procedure up to the final result $Du \in W^{\alpha \gamma,2}$ for every $\gamma < 1$.


\section{Structure conditions and notation}
\label{struc-not}

We impose on the coefficients $a \colon \Omega \times \R^N \times \R^{nN} \to \R^{nN}$ standard conditions (here stated for general $p$-growth, even if we will concentrate on case $p=2$): the mapping $z \mapsto a(x,u,z)$ is a continuous vector field, and for fixed numbers $0 < \nu \leq L$, $p \in (1,\infty)$ and all $x, \bar{x} \in \Omega$, $u,\bar{u} \in \R^N$, $z \in \R^{nN}$, the following growth, ellipticity and continuity assumptions hold:
\begin{align*}
\text{{\bf (H1)}} \qquad  & a \text{ has polynomial growth and is differentiable in } z \text{ with continuous, bounded derivatives:} \\
	& \hspace{1cm} \kabs{a(x,u,z)} +  \big(1 + \kabs{z}^2\big)^{\frac{1}{2}} 
	\, \kabs{ D_z a (x,u,z) }\, \leq \, L \, \big(1 + 
	\kabs{z}^{2}\big)^{\frac{p-1}{2}} \,, \\
\text{{\bf (H2)}} \qquad & a \text{ is uniformly strongly elliptic, i.\,e.} \\
	& \hspace{1cm} D_z a (x,u,z) \, \lambda \cdot \lambda \, \geq 
	\, \nu \, \big(1 + \kabs{z}^2\big)^{\frac{p-2}{2}} \kabs{\lambda}^2 \qquad \forall \,
	\lambda  \in \R^{nN} \,,\\
\text{{\bf (H3)}} \qquad & \text{There exists a nondecreasing, concave modulus of continuity } 
	\omega_{\alpha}:\R^+ \to [0,1] \\[-0.1cm]
	& \text{such that } \omega_{\alpha}(s) \leq \min\{1,s^{\alpha}\} \text{ for all } s \in \R^+ \text{ and}\\
	& \hspace{1cm} \kabs{a(x,u,z) -  a(\bar{x}, \bar{u}, z)} \, \leq \, L 
	\, \big(1 + \kabs{z}^2\big)^{\frac{p-1}{2}} \, \omega_{\alpha}\big(\kabs{x-\bar{x}} 
	+ \kabs{u-\bar{u}}\big) \,. \\
\intertext{The latter condition (H3) prescribes uniform H\"older continuity with respect to the $(x,u)$-variable with H\"older exponent $\alpha$ (for fixed $z$). Moreover, we assume the inhomogeneity $b \colon \Omega \times \R^N \times \R^{nN} \to \R^N$ to be a Carath\'eodory map, that is, it is continuous with respect to $(u,z)$ and measurable with respect to $x$, and to satisfy a natural growth condition of the form}
\text{{\bf (B)}} \qquad & \text{there exists a constant } L_2 
	\text{ (possibly depending on } M > 0 \text{) such that}  \\
	& \hspace{1cm} \kabs{b(x,u,z)} \leq L + L_2 \, \kabs{z}^p \\
	& \text{for all } x \in \Omega, u \in \R^N \text{ with } \kabs{u} \leq M, \text{ and } z \in \R^{nN} \,.
\end{align*}

We further make some remarks on the notation used below: 

\emph{(Half-)Balls, cubes and cylinders:} We write $B_{\rho}(y) = \{ x \in \R^n: \kabs{x-y}< \rho \}$ and $B_{\rho}^+(y) = \{ x \in \R^n: x_n > 0, \, \kabs{x-y} < \rho \}$ for an $n$-dimensional ball or the intersection of the ball with the upper half-space $\R^{n-1} \times \R^+$, centered at a point $y \in \R^n$ (respectively $\in \R^{n-1} \times \R^+_0$ in the latter case) with radius $\rho>0$. In the case $y = 0$ we set $B_{\rho}:= B_{\rho}(0)$, $B:=B_1$ as well as $B^+_{\rho}:= B^+_{\rho}(0)$, $B^+:= B^+_{1}$. Furthermore, we denote by $D_{\rho}(y')$ the $(n-1)$-dimensional ball $D_{\rho}(y') := \big\{ x \in \R^{n-1} \colon \kabs{y' - x'} < \rho \big\}$ for $y' \in \R^{n-1}$, and by $Z_{\rho}(y)$ the open cylinder on the upper half-plane $\R^{n-1} \times \R^+$
\begin{equation*}
Z_{\rho}(y) \, := \, D_{\rho}(y') \times \big(\max\{0,y_n-\rho\},y_n+\rho \big) 
	\, =: \,  D_{\rho}(y') \times I_{\rho}(y_n) 
\end{equation*}
for a center $y=:(y',y_n) \in \R^n$ with $y_n \geq 0$. Similarly as for balls, cubes with center $y \in \R^n$ and side-length $2\rho$ are denoted by $Q_{\rho}(y)$, upper half-cubes by $Q^+_{\rho}(y)$, and we further write $Q^0_{\rho}(y) = \partial Q^+_{\rho}(y) \cap \R^{n-1} \times \{0\}$ (with the corresponding abbreviations for $y=0$ and if $\rho =1$).

\emph{Function spaces:} We will work with functions belonging to the H\"older spaces $C^{1,\alpha}$, $\alpha \in (0,1)$, and the (fractional) Sobolev space $W^{\theta,p}$, $\theta \in (0,1], p \in [1,\infty)$. The definition for noninteger values of $\theta$ and some preliminary material is collected in the next section. Moreover, we introduce the following notation for $W^{1,p}$-functions defined on a upper half-cubes $Q^+_\rho(y)$ which vanish on $Q^0_{\rho}(y)$ (in the sense of traces):
\begin{equation*}
W^{1,p}_{\Gamma}(Q^+_\rho(y),\R^N) := \big\{ u \in W^{1,p}(Q^+_\rho(y),\R^N): u = 0 \text{ on } Q^0_{\rho}(y) \big\} \,.
\end{equation*}
where $y_n < \rho$ is satisfied; the subspace of functions vanishing on the whole boundary is denoted by $W^{1,p}_0$. Sometimes, it will be convenient to treat the tangential derivative $D'u := (D_1 u,\ldots,D_{n-1}u)$ and the normal derivative $D_n u$ of a Sobolev function $u$ separately.

\emph{Measures and mean values:} For a given set $X \subset \R^k$ we write $\Lk(X) = \kabs{X}$ and $\dim_{\mathcal{H}}(X)$ for its $k$-dimensional Lebesgue-measure and its Hausdorff dimension, respectively. Furthermore, if $h \in L^1(X, \R^N)$ and $0<\kabs{X} < \infty$, we denote the average of $h$ by $(h)_{X}=\tmint_X h \dx$, and when working on cylinders we will use the abbreviation $(v)_{x_0,\rho}:=(v)_{Z_{\rho}(x_0)}$. We further define the slice-wise mean value of $u$ in $D_{r}((x_0)')$ at almost every height $x_n \in I_\rho((x_0)_n)$ via
\begin{equation*}
(v)_{x_0',\rho}(x_n) \, := \, \mI{D_{\rho}((x_0)')} v(x',x_n) \dx'\,.
\end{equation*}

The constants $c$ appearing in the different estimates will all be chosen greater than or equal to $1$, and they may vary from line to line. 


\section{Fractional Sobolev spaces and finite differences}
\label{fractSob}

In what follows, we will use the notation of \cite{ADAMS75} (see also \cite{KRIMIN05,DUZKRIMIN05}). For a bounded open set $A \subset \R^n$, parameters $\theta \in (0,1)$ and $q \in [1,\infty)$ we write $u \in W^{\theta,q}(A,\R^N)$ provided that $u \in L^q(A,\R^N)$ and the following Gagliardo-type norm of $u$ defined as 
\begin{equation*}
\knorm{u}_{W^{\theta,q}(A)} \, := \, \Big( \int_A \kabs{u(x)}^q \dx \Big)^{\frac{1}{q}} + \Big( \int_A \int_A \frac{\kabs{u(x) - u(y)}^q}{\kabs{x-y}^{n + q \theta}} \dx \, dy \Big)^{\frac{1}{q}}
\end{equation*}
is finite. In order to formulate a general criterion for a function to belong to a fractional Sobolev space we introduce the finite difference operator $\tau_{e,h}$ with respect to a direction $e \in B_1 \subset \R^n$ and with stepsize $h \in \R$ via
\begin{equation*}
\tau_{e,h} G(x) \, \equiv \, \tau_{e,h} (G)(x) \, := \, G(x + h e) - G(x)
\end{equation*}
for a vector valued function $G: A \to \R^{N}$ (this makes sense whenever $x, x + h e \in A$). If $e = e_s$, $s \in \{1,\ldots,n \}$, is a standard basis vector, we use the abbreviation $\tau_{s,h}$ instead of $\tau_{e_s,h}$. These finite differences are related to the fractional Sobolev spaces (in the interior as well as in an up-to-the-boundary version) via the next lemma:

\begin{lemma}[\cite{KRIMIN05a}, Lemma 2.5; \cite{DUZKRIMIN05}, Lemma 2.2]
\label{frac-Sobol-lemma}
Let $G \! \in \! L^q(Q_R^{+},\R^{N})$, $q \geq 1$, and assume that for $\theta \in (0,1]$, $M > 0$ and some $0 < r < R$ we have
\begin{equation*}
\sum_{s=1}^n \int_{Q_r^{+}} \kabs{\tau_{s,h} G}^q \dx \, \leq \, M^q \, \kabs{h}^{q \theta}
\end{equation*}
for every $h \in \R$ satisfying $0 < \kabs{h} \leq d$ where $0 < d < \min\{1,R-r\}$ is a fixed number. In the case $s=n$ we only allow positive values of $h$. Then $G \in W^{b,q}(Q_{\rho}^{+},\R^{N})$ for every $b \in (0,\theta)$ and $\rho < r$. Moreover, there exists a constant $c = c(n,q)$ (in particular, independent of $M$ and $G$) such that the following inequality holds true:
\begin{equation*}
\int_{Q_{\rho}^{+}} \int_{Q_{\rho}^{+}} \frac{\kabs{G(x)-G(y)}^q}{\kabs{x-y}^{n+bq}} \dx \, dy \, \leq \, 
c \, \Big( \frac{M^q \epsilon^{q(\theta -b)}}{\theta - b} + \frac{\kabs{Q_R^{+}}}{\epsilon^{n + b q}} \int_{Q_{R}^{+}} \kabs{G}^q \dx \Big) \,,
\end{equation*}
where $\epsilon := \min \{r - \rho,d \}$. In the interior the same result holds true without any constraint on the sign of $h$ with respect to the direction of the differences $\tau_{s,h}$. Moreover, we can consider (half-)balls or cylinders instead of cubes.
\end{lemma}

In the case where $G$ is the weak derivative of a $W^{1,q}$ function $v$ and where an estimate for finite differences only in tangential direction is known, we are still in a position to state a fractional differentiability result which is limited to the tangential derivative of $v$:

\begin{lemma}
\label{frac-Sobol-lemma-2}
Let $v \in W^{1,q}(Q_R^{+},\R^{N})$, $q \geq 1$, and assume that for $\theta \in (0,1]$, $M > 0$ and some $0 < r < R$ we have
\begin{equation}
\label{assumption-frac-Sobol-lemma-2}
\sum_{s=1}^{n-1} \int_{Q_r^{+}} \kabs{\tau_{s,h} Dv}^q \dx \, \leq \, M^q \, \kabs{h}^{q \theta}
\end{equation}
for every $h \in \R$ satisfying $0 < \kabs{h} \leq d$ where $0 < d < \min\{1,R-r\}$ is a fixed number. Then $D'v = (D_1 v,\ldots,D_{n-1} v) \in W^{b,q}(Q_{\rho}^{+},\R^{(n-1)N})$ for every $b \in (0,\theta)$ and $\rho < r$.

\begin{proof}
We first fix $b \in (0,\theta)$ and $\rho \in (0, r)$. We consider arbitrary numbers $h' \in \R^+$ and $h \in \R$ satisfying $0 < \kabs{h},\kabs{h'} < \min\{d,\frac{r-\rho}{3}\}$. Then, using Young's inequality, standard properties of the difference operator and the assumption~(\ref{assumption-frac-Sobol-lemma-2}) on finite differences in tangential direction, we conclude for every $\epsilon \in (0,\theta)$ and $s \in \{1,\ldots,n-1\}$:
\begin{align*}
\lefteqn{ \hspace{-0.75cm} \kabs{h'}^{- (\theta - \epsilon) q} \, \kabs{h}^{-(1+\epsilon)q} 
	\int_{Q_{r-2d}^{+}} \kabs{\tau_{n,h'} \tau_{s,h} \tau_{s,-h} v}^q \dx } \\
	& \leq \, \big( \kabs{h'}^{-q} \,\kabs{h}^{-\theta q} + \kabs{h}^{-q-\theta q} \big)
	\int_{Q_{r-2d}^{+}} \kabs{\tau_{n,h'} \tau_{s,h} \tau_{s,-h} v}^q \dx \\
	& \leq \, 2 \, \kabs{h'}^{-q} \, \kabs{h}^{-\theta q} 
	\int_{Q_{r-d}^{+}} \kabs{\tau_{s,h} \tau_{n,h'} v}^q \dx
	+ 2 \, \kabs{h}^{-q-\theta q}
	\int_{Q_{r-d}^{+}} \kabs{\tau_{s,h} \tau_{s,-h} v}^q \dx \\
	& \leq \,  2 \, \kabs{h}^{-\theta q} 
	\int_{Q_{r}^{+}} \kabs{\tau_{s,h} D_n v}^q \dx
	+ 2 \, \kabs{h}^{-\theta q}
	\int_{Q_{r}^{+}} \kabs{\tau_{s,h} D_s v}^q \dx \, \leq \, 4 \, M^q
\end{align*}
uniformly in $h$, $h'$. From \cite[Lemma 2.2.1]{DOMOKOS04} we infer (for possibly smaller values of $\kabs{h}$)
\begin{equation*}
\kabs{h'}^{- (\theta - \epsilon) q} \, \kabs{h}^{-q} 
	\int_{Q_{r-2d}^{+}} \kabs{\tau_{n,h'} \tau_{s,h} v}^q \dx
	\, \leq \, c \, \Big( \int_{Q_R^+} \kabs{Dv}^q \dx + M^q \Big) \,,
\end{equation*}
and the constant $c$ depends only on $\theta, q, \epsilon, d$ and $r-\rho$. Considering the limit $h \to 0$, we hence end up with
\begin{equation*}
\kabs{h'}^{- (\theta - \epsilon) q} 
	\int_{Q_{r-2d}^{+}} \kabs{\tau_{n,h'} D_s v}^q \dx
	\, \leq \, c \, \Big( \int_{Q_R^+} \kabs{Dv}^q \dx + M^q \Big) \,.
\end{equation*}
Keeping in mind that the index $s \in \{1,\ldots,n-1\}$ is arbitrary, we combine the latter inequality with \eqref{assumption-frac-Sobol-lemma-2} to find
\begin{equation*}
\sum_{s=1}^{n} \int_{Q_{r-2d}^{+}} \kabs{\tau_{s,h} D'v}^q \dx \, \leq \, 
	c \, \kabs{h}^{(\theta - \epsilon) q} \, \Big( \int_{Q_R^+} \kabs{Dv}^q \dx + M^q \Big)
\end{equation*}
for all $h \in \R$ satisfying $0 < \kabs{h} \leq \min\{d,\frac{r-\rho}{3}\}$ where we only allow positive values of $h$ if $s=n$. For $\epsilon = (\theta - b) /2$ the application of Lemma~\ref{frac-Sobol-lemma} with $\theta,r$ replaced by $\theta - \epsilon, r-2d$ finishes the proof.
\end{proof}
\end{lemma}

The following interpolation inequality can be found in \cite[Lemma 2.V]{CAMPANATO82a} and is essentially based on the inequality in \cite[Theorem 2.I]{CAMPCAN81} for the case $p=2$.

\begin{theorem}
\label{CampInterpolation}
Let $\lambda, \theta \in (0,1]$, $p \in (1,\infty)$ and $u \in C^{0,\lambda}(\overline{Q},\R^N)$ such that $Du \in W^{\theta,p}(Q,\R^{nN})$ with $p \theta < n$, where $Q \subset \R^N$ is an (upper) cube. Then
\begin{equation*}
Du \in L^s(Q,\R^{nN}) \hspace{2cm} \text{ for all } s < \frac{np(1+\theta)}{n-p \theta \lambda} \,.
\end{equation*}
Moreover,
\begin{equation*}
\int_Q \kabs{Du}^s \dx \, \leq \, c\big( n,N,p,\theta,\lambda,s,\kabs{Q}, \knorm{u}_{W^{1+\theta,p}(Q,\R^N)},[u]_{C^{0,\lambda}(\overline{Q},\R^N)}\big) \,.
\end{equation*}
\end{theorem}

The next lemma enables us to conclude from difference estimates for a map $v$ an appropriate estimate for the averaged mean deviation with respect to slice-wise mean values:

\begin{lemma}[\cite{KRONZHABIL}]
\label{rbp-manni-slice-lemma}
Let $\sigma < \frac{1}{3}$, $n \geq 2$, $\tau>0$, $Z_{\rho}(x_0) \subset Q^+$ for some $x_0 \in Q^+ \cup Q^0$. Furthermore, assume that $v \in L^p(Z_{\rho}(x_0),\R^{N})$, $p>1$, satisfies
\begin{equation*}
\int_{Z_{\sigma \rho}(x_0)} \kabs{\tau_{h,e} v}^p \dx \, \leq \, K^p \, \kabs{h}^{\tau p}
\end{equation*}
for some $K>0$, all $e \in S^{n-1}$ with $e \perp e_n$ and $h \in \R$ with $\kabs{h} < 2 \sigma \rho$. Then, for every $\beta \in (0,\tau)$ there exists a function $F \in L^p(Z_{\sigma \rho}(x_0))$ such that
\begin{equation*}
\int_{Z_{\sigma \rho}(x_0)} \kabs{F}^p \dx \, \leq \, c(n,p,\tau,\beta) \, K^p \, \rho^{(\tau-\beta)p}
\end{equation*}
and
\begin{equation*}
\Big( \mI{Z_{r}(z)} \mI{D_{r}(z')} \kabs{v(x',x_n) - v(y',x_n)}^ {\widetilde{p}}
	\dy' \dx \Big)^{\frac{1}{\widetilde{p}}} \, \leq \, c(n,\beta) \, 
	r^{\beta} \, F(z)
\end{equation*}
for every exponent $\widetilde{p} \in [1,p)$, almost all $z \in Q^+ \cup Q^0$ and all $r>0$ such that $Z_{r}(z) \subset Z_{\sigma \rho}(x_0)$.
\end{lemma}

A different definition for fractional Sobolev spaces, based on pointwise inequalities, can be derived as follows: Let $\Omega \subset \R^n$ be a bounded domain, $p \geq 1$ and $\theta \in (0,1]$. Following the approach of Haj{\l}asz in \cite{HAJLASZ96}, we set
\begin{align*}
{\cal D}^{\theta,p}(\Omega;f) & := \, \big\{ g \in L^p(\Omega) \colon \exists \, 
	E \subset \Omega, \, \kabs{E} = 0 \text{ such that} \\
	& \hspace{2.9cm} \kabs{f(x)-f(y)} \leq \kabs{x-y}^{\theta} (g(x)+g(y)) \text{ for all } 
	x,y \in \Omega \setminus E \big\}\,,
\end{align*}
and we define the fractional Sobolev space via
\begin{equation*}
M^{\theta,p}(\Omega,\R^N) \, := \, \big\{ f \in L^p(\Omega,\R^N) \colon {\cal D}^{\theta,p}(\Omega;f) \neq \emptyset \big\} \,.
\end{equation*}
$M^{\theta,p}(\Omega,\R^N)$ is equipped with the norm 
\begin{equation*}
\knorm{f}_{M^{\theta,p}(\Omega,\R^N)} \, := \, \knorm{f}_{L^p(\Omega,\R^N)} + \inf_{g \in {\cal D}^{\theta,p}(\Omega;f)} \knorm{g}_{L^p(\Omega)} \,.
\end{equation*}
For $p \in (1,\infty)$, due to the convexity of $L^p$, to every $f \in M^{\theta,p}(\Omega,\R^N)$ there exists a unique function $g \in L^p(\Omega)$ which minimizes the $L^p(\Omega)$-norm amongst all functions in ${\cal D}^{\theta,p}(\Omega;f)$. We highlight that this definition has its origin in the definition of Sobolev spaces in the context of arbitrary metric spaces (replacing $\kabs{x-y}$ by $\dist(x,y)$) and that it does not use of the notion of derivatives (for a more detailed discussion of the metric setting we refer to \cite{HAJKOS00}). Employing the Hardy-Littlewood maximal function we see that this ``metric'' Sobolev space coincides with the classical Sobolev space for the integer order $\theta=1$ and sufficiently regular domains (e.\,g. with Lipschitz boundary). More precisely, provided that $p>1$, there holds $M^{1,p}(\Omega,\R^N) \, = \, W^{1,p}(\Omega,\R^N)$
for all bounded domains $\Omega$ with the so-called extension property, meaning that there exists a bounded linear operator $E: W^{1,p}(\Omega,\R^N) \to W^{1,p}(\R^n,\R^N)$ such that for every $f \in W^{1,p}(\Omega,\R^N)$ there holds $Ef=f$ almost everywhere in $\Omega$. Instead, the equivalence fails if $p=1$, see \cite{HAJLASZ95}. Furthermore, the definitions of the classical and the metric fractional Sobolev spaces immediately yield for all bounded domains $\Omega$, fractional orders $\theta \in (0,1)$ and $p \in [1,\infty)$ the following inclusion:
\begin{equation*}
M^{\theta,p}(\Omega,\R^N) \subseteq W^{\theta^{\prime},p}(\Omega,\R^N) \qquad \text{ for all } \theta^{\prime} \in (0,\theta).
\end{equation*}

The following lemma provides an integral characterization of fractional Sobolev spaces for domains satisfying the mild Ahlfors regularity condition, which demands the existence of a positive constant $k_{\Omega}$ such that
\begin{equation*}
(K_{\Omega}) \qquad \kabs{B_{\rho}(x_0) \cap \Omega} \geq k_{\Omega} \, \rho^n \qquad \text{ for all points } x_0 \in \overline{\Omega} \text{ and every radius } \rho \leq \diam(\Omega)\,.
\end{equation*}
In other words: the domain is not allowed to have external cusps. We note that the latter condition is for example satisfied by the large class of domains with Lipschitz-continuous boundary.

\begin{lemma}
\label{lemma-frac-char}
Let $\Omega \subset \R^n$ be a domain which fulfills an Ahlfors condition $(K_{\Omega})$, $\theta \in (0,1]$, $p \in (1,\infty)$. Then the following two statements are equivalent:
\begin{itemize}
\item[(i)] $f \in M^{\theta,p}(\Omega,\R^N)$
\item[(ii)] $f \in L^1(\Omega,\R^N)$ and there exists a function $h \in L^p(\Omega)$ and a radius $R_0 > 0$ such that
\begin{equation}
\label{frac-equiv-lemma}
\mI{B_{\rho}(x_0) \cap \Omega} \kabs{f - (f)_{B_{\rho}(x_0) \cap \Omega}} \dx \, \leq \, \rho^{\theta} \, h(x_0)
\end{equation}
for almost all $x_0 \in \overline{\Omega}$ and $\rho \leq R_0$.
\end{itemize}

\begin{proof}
The implication (i) $\Rightarrow$ (ii) follows by standard properties of the Hardy-Littlewood maximal function for the choice $h = 4 M(g)$ with $g \in D^{\theta,p}(\Omega;f)$. The reverse implication (ii) $\Rightarrow$ (i) is an easy adaptation of the proof of Campanato's integral characterization of H\"older continuous functions, see e.\,g. \cite[Chapt. 1.1, Lemma 1]{SIMON96}.
\end{proof}
\end{lemma} 

\begin{remarks}
\label{remark-frac-char}
In fact, the following local version of the integral characterization holds: let $x_0 \in \overline{\Omega}$ and $R>0$ such that 
\begin{equation*}
\mI{B_r(z) \cap \Omega} \kabs{f - (f)_{B_{r}(z) \cap \Omega}} \dx \, \leq \, r^{\theta} \, h(z)
\end{equation*}
for almost all $z \in \overline{\Omega}$, $B_r(z) \subset B_{R}(x_0)$ and $h \in L^p(\Omega)$ as above. Then there holds $f \in M^{\theta,p}(B_{R/2}(x_0) \cap \Omega,\R^N)$ with
\begin{equation*}
\kabs{f(x)-f(y)} \, \leq \, c(n,k_{\Omega},\theta) \, \kabs{x-y}^{\theta} \, \big( h(x) + h(y)\big)
\end{equation*}
for almost all $x,y \in B_{R/2}(x_0) \cap \Omega$. In view of Jensen's inequality and the fact that the Hardy Littlewood maximal operator is a bounded map from $L^p$ to itself, this characterization allows to infer the inclusion
\begin{equation*}
W^{\theta,p}(\Omega,\R^N) \subseteq M^{\theta,p}(\Omega,\R^N)
\end{equation*}
whenever $\Omega$ satisfies an Ahlfors condition $(K_{\Omega})$, $\theta \in (0,1)$ and $p \in (1,\infty)$.

Moreover, we note that (i) implies indeed the following statement: there exists a function $h \in L^p(\Omega)$ and a radius $R_0 > 0$ such that
\begin{equation*}
\Big(\mI{B_{\rho}(x_0) \cap \Omega} \kabs{f - (f)_{B_{\rho}(x_0) \cap \Omega}}^q \dx\Big)^{\frac{1}{q}} \, \leq \, \rho^{\theta} \, h(x_0)
\end{equation*}
for all $q<p$ and almost all $x_0 \in \overline{\Omega}$ and $\rho \leq R_0$.
\end{remarks}


\section{Some basic facts about the solution}

In what follows, we restrict ourselves to the model case $\Omega=Q_2^+$, and we study weak solutions $u \in W^{1,p}_{\Gamma}(Q_2^+,\R^N) \cap L^{\infty}(Q_2^+,\R^N)$ of the system
\begin{equation}
\label{rbp-manni-model}
- \diverg a(\ccdot,u,Du) \, = \, b(\ccdot,u,Du) \qquad \text{in } Q_2^+. 
\end{equation}
By a transformation argument this covers the situation of general inhomogeneous systems of type \eqref{DP-rbp} on arbitrary domains $\Omega$ of class $C^{1,\alpha}$. Moreover, we argue under the permanent assumption that the weak solution $u$ of system \eqref{rbp-manni-model} is H\"older continuous  on $Q^+$ with H\"older exponent $\lambda$ for some $\lambda \in (0,1)$. This assumption will later be justified by the fact that in low dimensions the weak solution $u$ is a~priori known to be H\"older continuous outside a set of Hausdorff dimension $n-2$ (and since we are interested in the behavior of $Du$ on the boundary which is of Hausdorff dimension $n-1$ this information is sufficient to forget about the bad set where $u$ is not H\"older continuous). 

We now present some tools needed in the remainder of the paper: first, we recall the well-known Caccioppoli inequality in an up-to-the-boundary version. The fact that the oscillations of $u$ are due to its continuity arbitrarily small in a cylinder -- provided that the side length of the cylinder is chosen sufficiently small~-- allows to simplify the estimates which are usually slightly more involved for nonlinear elliptic systems with inhomogeneities under a natural growth condition. As a matter of fact we here do not need the smallness assumption $\kabs{u} \leq M$ with $2 L_2 M  < \nu$. 

\begin{lemma}[Caccioppoli inequality revised]
\label{manni-caccio-rev}
Let $u \in W^{1,p}_{\Gamma}(Q_2^+,\R^N) \cap L^{\infty}(Q_2^+,\R^N)$ be a weak solution of \eqref{rbp-manni-model} under the assumptions (H1)-(H3) and (B). Assume further $u \in C^{0,\lambda}(Q^+,\R^N)$. Then there exist positive constants $\widetilde{c}_{cacc} = \widetilde{c}_{cacc}(n,N,p,\tfrac{L}{\nu},\tfrac{L_2}{\nu})$ and $\widetilde{\rho}_{cacc} =\widetilde{\rho}_{cacc}(p,\tfrac{L}{\nu},\tfrac{L_2}{\nu}, \lambda,[u]_{C^{0,\lambda}(Q^+,\R^N)})$ such that for every $\xi \in \R^N$ and every cylinder $Z_{\rho}(y) \subset Q^+$ with $y \in Q^+ \cup Q^0$ and $y_n < \rho \leq \widetilde{\rho}_{cacc}$ there holds:
\begin{equation*}
\mI{Z_{\rho/2}(y)} \kabs{V(Du)-V(\xi \otimes e_n)}^2 \dx \\
	 \leq \, \widetilde{c}_{cacc} \, \Big( \mI{Z_{\rho}(y)} \Babs{ 
	V\Big(\frac{u - \xi  x_n}{\rho} \Big)}^2 \dx + 
	\rho^{2 \alpha} \, \big( 1 + \kabs{\xi} \big)^{p+2\alpha} \Big) \,.
\end{equation*}
\end{lemma}

Here we have used the $V$-function which is in general defined by $V(\xi) = (1 + \kabs{\xi}^2)^{(p-2)/4} \xi$ for all $\xi \in \R^k$ for some $k \in  \N$ (in the quadratic case it is just the identity map) and which is in particular a bi-Lipschitz bijection on $\R^k$. Secondly, we recall an estimate concerning finite tangential differences of $Du$ which is the starting point to proceed to fractional differentiability estimates for $Du$ and hence to dimension reduction arguments for the singular set: We consider $\delta \in (0,1)$ and assume  $u \in W^{1,p}_{\Gamma}(Q_2^+,\R^N) \cap L^{\infty}(Q_2^+,\R^N)$ to be a weak solution of system \eqref{rbp-manni-model}. Then for every cut-off function $\eta \in C_0^{\infty}(Q_{1-\delta},[0,1])$ and every tangential direction $e \in S^{n-1}$ with $e \perp e_n$ there holds
\begin{align}
\label{rbp-k-start}
\int_{Q^+} \eta^2 \kabs{\tau_{e,h} V(Du)}^2 \dx & \leq \,  c \, \bigg( \kabs{h}^{2 \alpha} \int_{Q^+ \cap \supp(\eta)} 
	\big( 1 + \kabs{Du(x)}^p + \kabs{Du(x+h e)}^p +  \kabs{h}^{-p} \kabs{\tau_{e,h}u(x)}^p \big) 
	\dx \nonumber \\
	& \qquad {}+ \int_{Q^+ \cap \supp(\eta)}  
	\big( 1 + \kabs{Du(x)}^2 + \kabs{Du(x + h e_s)}^2 \big)^{\frac{p}{2}} 
	\, \kabs{\tau_{e,h} u(x)}^{2 \alpha} \dx \nonumber \\
	& \qquad {}+  \int_{Q^+} \big( 1 + \kabs{Du(x)}^p\big) 
	\, \kabs{\tau_{e,-h}(\eta^2 \tau_{e,h} u(x))} \dx \bigg)
\end{align}
for all $h \in \R$ with $\kabs{h} < \delta$, and the constant $c$ depends only on $n,N,p,\tfrac{L}{\nu},\tfrac{L_2}{\nu}$, $\knorm{u}_{L^\infty}$ and $\knorm{D\eta}_{L^\infty}$. We highlight that this estimate is the up to the boundary analogue of \cite[estimate (4.7)]{MINGIONE03}, and its proof follows the line of arguments in \cite{MINGIONE03}: Testing the weak formulation of \eqref{rbp-manni-model} with the function $\tau_{e,-h} \varphi$ for $\varphi \in W^{1,p}_0(Q^+,\R^N) \cap L^{\infty}(Q^+,\R^N)$ with $\supp \varphi \subset Q_{1-\delta}$, we first use partial integration for finite differences on the left-hand side which results in integrals involving $\tau_{s,h} \big(a(x,u(x),Du(x))\big)$. Decomposing
\begin{align}
\label{def-ABC}
\lefteqn{\hspace{-0.75cm} \tau_{e,h} \big(a(x,u(x),Du(x))\big)} \nonumber  \\
& = \,  a(x+h e,u(x+h e),Du(x+h e)) -  a(x,u(x+h e),Du(x+h e)) \nonumber \\
& \qquad {} + a(x,u(x+h e),Du(x+h e)) -  a(x,u(x),Du(x+h e)) \nonumber \\
& \qquad {} +  a(x,u(x),Du(x+h e)) -  a(x,u(x),Du(x)) \nonumber \\
& =: \, {\cal A} (h) + {\cal B} (h) + {\cal C} (h) \,,
\end{align}
we hence find
\begin{equation}
\label{test-varphi}
\int_{Q^+} \big[ {\cal A} (h) + {\cal B} (h) + {\cal C} (h)\big] \cdot D\varphi \dx  \, 
= \, \int_{Q^+} b(x,u,Du) \cdot \tau_{e,-h} \varphi \dx \,.
\end{equation}
Choosing $\varphi = \eta^2 \tau_{e,h} u$, we have to estimate the various terms by taking advantage of the growth and continuity assumptions of the coefficients and the inhomogeneity exactly as in \cite{MINGIONE03}, and we then end up with the desired inequality \eqref{rbp-k-start}.


\section{The proof of Theorem~\ref{regul-bp-mit-u-nat-2}}

\subsection{Higher integrability of finite differences}
\label{high-anfang}

We first state a higher integrability estimate for both $Du$ and for finite differences of $Du$ (again motivated from \cite{MINGIONE03}), which will allow later to end up with a slightly sharper estimate on the Hausdorff dimension of the singular set. We first observe the well-known existence of a higher integrability exponent $s_0>2$ depending only on $n,N,\frac{L}{\nu}, \frac{L_2}{\nu}$ and $[u]_{C^{0,\lambda}(Q^+,\R^N)}$ such that $u \in W^{1,s_0}(Q^+_{\rho},\R^N)$ for all $\rho<1$.  Furthermore, for every center $x_0 \in Q^+ \cup Q^0$ and every radius  $\rho \in (0,1-\kabs{x_0})$ there holds 
\begin{align}
\label{rbp-manni-higher-anfang}
\Big( \mI{Z_{\rho/2}(x_0)} \kabs{Du}^{s_0} \dx \Big)^{\frac{1}{s_0}}
\, \leq \, c\big(n,N,\tfrac{L}{\nu},\tfrac{L_2}{\nu},[u]_{C^{0,\lambda}(Q^+,\R^N)}\big) \, \Big( \mI{Z_{\rho}(x_0)} \big( 1 + \kabs{Du}^2\big) \dx \Big)^{\frac{1}{2}} \,,
\end{align}
see e.g. \cite[Lemma 4.1]{BECK09b}. Combining the higher integrability with \eqref{rbp-k-start} we obtain similarly to \cite[Section 5, step 2]{MINGIONE03} a higher integrability result for $\tau_{e,h}Du$:

\begin{proposition}
\label{rbp-manni-diff-higher-anfang}
Let $u \in W^{1,2}_{\Gamma}(Q_2^+,\R^N) \cap L^{\infty}(Q_2^+,\R^N) \cap C^{0,\lambda}(Q^+,\R^N)$ be a weak solution of \eqref{rbp-manni-model} under the assumptions (H1)-(H3) and (B). Furthermore, let $Z_{\rho}(x_0) \subset Q^+$ for some $x_0 \in Q^+ \cup Q^0$, $\sigma \in (0,\frac{1}{10})$, $e \in S^{n-1}$ with $e \perp e_n$ and $h \in \R$ with $\kabs{h} \in (0,2\sigma\rho)$. Then there exists a higher integrability exponent $s \in (2,s_0)$ depending only on $n,N,\frac{L}{\nu},\frac{L_2}{\nu}$ and $[u]_{C^{0,\lambda}(Q^+,\R^N)}$ such that
\begin{equation*}
\mI{Z_{\sigma \rho}(x_0)} \kabs{\tau_{e,h}Du}^{s} \dx 
	\, \leq \, c \, \kabs{h}^{\frac{\alpha \lambda s}{2}} \, 
	\Big(\mI{Z_{\rho}(x_0)} \big( 1 + \kabs{Du}^2 
	\big) \dx \Big)^{\frac{s}{2}}
\end{equation*}
for a constant $c \!= c\big(n,N,\frac{L}{\nu},\frac{L_2}{\nu},[u]_{C^{0,\lambda}(Q^+,\R^N)}, \rho,\sigma\big)$.

\begin{proof} 
We consider in the sequel the tangential directions $e \in S^{n-1}$, i.\,e. $e \perp e_n$, and we initially look at numbers $h \in \R$ satisfying $\kabs{h}<1$. Recalling the abbreviations for ${\cal A} (h)$, ${\cal B} (h)$ and ${\cal C} (h)$ from \eqref{def-ABC}, representing the differences of the coefficients $a(\cdot,\cdot,\cdot)$ with respect to each variable, we set
\begin{equation*}
v_h \, := \, \frac{\tau_{e,h}u}{\kabs{h}^{\frac{\alpha \lambda}{2}}} \,, \hspace{1cm}
\widetilde{{\cal A}}(h) \, := \, \frac{-{\cal A} (h)}{\kabs{h}^{\frac{\alpha \lambda}{2}}} \,, \hspace{1cm}
\widetilde{{\cal B}}(h) \, := \, \frac{-{\cal B} (h)}{\kabs{h}^{\frac{\alpha \lambda}{2}}} \,,
\end{equation*}
and we define $\widetilde{{\cal C}}(h) = \int_0^1 D_z a\big(x,u(x),Du(x) + t\tau_{e,h} Du(x))\big) \, dt$. Dividing the previous identity \eqref{test-varphi} by $\kabs{h}^{\alpha \lambda/2}$ (which is half the power of $\kabs{h}$ to be expected in \eqref{rbp-k-start} for $\lambda$-H\"older continuous solutions) we get
\begin{equation}
\label{system-v-h}
\int_{Q^+} \widetilde{{\cal C}}(h) \, Dv_h \cdot D\varphi \dx  \, 
= \, \int_{Q^+} \big[\widetilde{{\cal A}}(h) + \widetilde{{\cal B}}(h) \big] \cdot D\varphi \dx + 
	\int_{Q^+} \kabs{h}^{-\frac{\alpha \lambda}{2}} \, b(x,u,Du) \cdot \tau_{e,-h} \varphi \dx
\end{equation}
for all functions $\varphi \in W^{1,2}_0(Q^+_{1-\kabs{h}},\R^N) \cap L^{\infty}(Q^+_{1-\kabs{h}},\R^N)$, i.\,e. the map $v_h \in W^{1,2}(Q^+_{1-\kabs{h}},\R^N)$ is a weak solution to the linear system \eqref{system-v-h} for every $h \in \R$ with $\kabs{h}<1$. In the next step we infer Caccioppoli-type inequalities for the functions $v_h$, for which the constants may be chosen independently of the parameter $h$. For this purpose we first observe some simple properties due to (H1)-(H3) and the H\"older continuity of $u$ with exponent $\lambda$:
\begin{align*}
\kabs{\widetilde{{\cal A}}(h)} & \leq \, L \, \big( 1 + \kabs{Du(x + h e)} \big) \,,\\
\kabs{\widetilde{{\cal B}}(h)} & \leq \, L \, [u]^{\alpha}_{C^{0,\lambda}(B^+,\R^N)}
	\, \big( 1 + \kabs{Du(x + h e)} \big) \,, \\
\nu \, \kabs{\widetilde{\lambda}}^2 & \leq \, \widetilde{{\cal C}}(h) \widetilde{\lambda} 
	\otimes \widetilde{\lambda} \, \leq \, L \, \kabs{\widetilde{\lambda}}^2  
	\qquad \forall \, \widetilde{\lambda} \in \R^{nN}\,.
\end{align*}
For $\sigma, \rho$ and $x_0$ fixed according to the assumptions of the proposition, we next choose $h \in \R$ such that $\kabs{h} \in (0,2\sigma\rho)$ and consider intersections of balls $B_R^+(y)$ with the upper half-plane $\R^{n-1} \times \R^+$ for centers $y \in \overline{Z_{(1-\sigma) \rho/2}(x_0)}$ satisfying $B_R^+(y) \subset Q_{1-\kabs{h}}^+$ (implying that $0 < R < 1-\kabs{h}-\max_{k \in \{1,\ldots,n\}}\kabs{y_k} $) and $y_n \leq \frac{3R}{4}$, i.\,e. we first study the situation for centers close to the boundary. Furthermore, we take a cut-off function $\eta \in C_0^{\infty}(B_{3R/4}(y),[0,1])$ satisfying $\eta \equiv 1$ on $B_{R/2}(y)$ and $\kabs{D\eta} \leq \frac{8}{R}$, and we choose $\varphi := \eta^2 v_h$ as a test function in \eqref{system-v-h}. Taking into account $D \varphi = \, \eta^2 Dv_h + 2 \eta v_h \otimes D\eta$, we estimate the various terms arising in \eqref{system-v-h}: using Young's inequality with $\epsilon \in (0,1)$ and the estimates for $\widetilde{{\cal A}} (h)$, $\widetilde{{\cal B}} (h)$ and $\widetilde{{\cal C}} (h)$ given above we see
\begin{itemize}
\item $\displaystyle \nu \int_{B_R^+(y)} \eta^2 \, \kabs{Dv_h}^2 \dx 
	\, \leq \, \int_{B_R^+(y)} \eta^2 \, \widetilde{{\cal C}} (h) \, Dv_h \cdot Dv_h \dx$\,,
\item $\displaystyle \int_{B_R^+(y)} 2 \, \eta \, \kabs{\widetilde{{\cal C}} (h) Dv_h \cdot v_h \otimes D\eta} \dx 
	\, \leq \, \epsilon \int_{B_R^+(y)} \eta^2 \, \kabs{Dv_h}^2 \dx 
	+ \frac{c \, L^2}{\epsilon \, R^2} \int_{B_R^+(y)} \kabs{v_h}^2 \dx$\,,
\item $\displaystyle \int_{B_R^+(y)} \kabs{\widetilde{{\cal A}} (h) \cdot D\varphi} \dx 
	\, \leq \, \epsilon \int_{B_R^+(y)} \eta^2 \, \kabs{Dv_h}^2 \dx 
	+ \frac{L}{R^2} \int_{B_R^+(y)} \kabs{v_h}^2 \dx$ \\[0.1cm]
	$\displaystyle \hspace{4.25cm} {}+ c \, \big(\epsilon^{-1} L^2 + L \big) \, 
	\int_{B_R^+(y)} \! \big( 1 + \kabs{Du(x + h e)}^2 \big) \dx$\,,
\item $\displaystyle \int_{B_R^+(y)} \kabs{\widetilde{{\cal B}} (h) \cdot D\varphi} \dx
	\, \leq \, \epsilon \int_{B_R^+(y)} \eta^2 \, \kabs{Dv_h}^2 \dx 
	+  \frac{c \, \epsilon}{R^2} \int_{B_R^+(y)} \kabs{v_h}^2 \dx$ \\[0.1cm]
	$\displaystyle \hspace{4.25cm} {}+ c\big([u]_{C^{0,\lambda}(B^+,\R^N)}\big) \, 
	\epsilon^{-1} L^2 \, 
	\int_{B_R^+(y)} \! \big( 1 + \kabs{Du(x + h e)}^2 \big) \dx$\,.
\end{itemize}
In order to estimate the last integral on the right-hand side of \eqref{system-v-h} we calculate
\begin{equation}
\label{doppelt-Absch}
\babs{\tau_{e,-h} \varphi} \,  = \, \babs{\tau_{e,-h}(\eta^2 v_h)} 
	\, \leq \, \kabs{h}^{-\frac{\alpha \lambda}{2}} 
	\big( \kabs{\tau_{e,h} u(x-he)} + \kabs{\tau_{e,h} u(x)}\big)
	\, \leq \, 2 \, [u]_{C^{0,\lambda}(Q^+,\R^N)} \, 
	\kabs{h}^{\lambda-\frac{\alpha \lambda}{2}} \,.
\end{equation}
This yields
\begin{itemize}
\item $\displaystyle \int_{B_R^+(y)} \kabs{h}^{-\frac{\alpha \lambda}{2}} \, \kabs{b(x,u,Du) \cdot \tau_{e,-h} \varphi} \dx \, \leq \, c\big([u]_{C^{0,\lambda}(Q^+,\R^N)}\big)
	\int_{B_R^+(y)} \big( L + L_2 \, \kabs{Du(x)}^2 \big) \dx$\,.
\end{itemize}
Collecting the estimates for all terms arising in equation \eqref{system-v-h} and choosing $\epsilon = \frac{\nu}{6}$, we finally conclude the Caccioppoli-type estimate
\begin{equation*}
\int_{B_{R/2}^+(y)}  \kabs{Dv_h}^2 \dx \, \leq \, c \, R^{-2} \int_{B_R^+(y)} \kabs{v_h}^2 \dx 
	+ c \int_{B_R^+(y)} \big( 1 + \kabs{Du(x)}^2 + \kabs{Du(x + h e)}^2 \big) \dx \,,
\end{equation*}
and the constant $c$ depends only on $\frac{L}{\nu},\frac{L_2}{\nu}$ and $[u]_{C^{0,\lambda}(Q^+,\R^N)}$. With the boundary version of the Sobolev-Poincar\'{e} inequality we deduce
\begin{align}
\label{v_h-bound}
\mI{B_{R/2}^+(y)} \kabs{Dv_h}^2 \dx 
	& \leq \, c \, \Big( \mI{B_R^+(y)} \kabs{Dv_h}^{\frac{2n}{n+2}} 
	\dx \Big)^{\frac{n+2}{n}} \nonumber \\
	& \qquad {}+ c \mI{B_R^+(y)} \big( 1 + 
	\kabs{Du(x)}^2 + \kabs{Du(x + h e)}^2 \big) \dx \,,
\end{align}
and the constant $c$ now depends additionally on the dimensions $n,N$. We here note that the integrand of the second integral on the right-hand side of the last inequality belongs to $L^{s_0/2}$ due to the higher integrability result for $Du$ from \eqref{rbp-manni-higher-anfang}. In the interior we proceed analogously and consider $B_R^+(y)$ with centers $y \in Z_{(1-\sigma) \rho/2}(x_0)$ satisfying $B_R^+(y) \subset Q_{1-\kabs{h}}^+$ and $y_n > \frac{3R}{4}$. If we choose $\varphi := \eta^2 \big(v_h - (v_h)_{y,3R/4} \big)$ as a test function all the computations above remain valid (with 2 replaced by 4 in inequality \eqref{doppelt-Absch}). Then, after applying the Sobolev-Poincar\'{e} inequality in the interior in the mean value version on the ball $B_{3R/4}(y)$, we obtain the corresponding inequality \eqref{v_h-bound} with the full ball $B_{R/2}(y)$ instead of $B_{R/2}^+(y)$, and $c$ has exactly the same dependencies as in the previous reverse H\"older-type inequality; in particular, the constant $c$ is independent of the parameter $h$.

Applying the global Gehring Lemma \cite[Theorem 2.4]{DUGROKRO04} on the cylinder $Z_{(1-\sigma) \rho/2}(x_0)$ for the choices of $\sigma, \rho$ and $x_0$ made in the assumptions of the proposition, we find that there exist a constant $c$ depending only on $n,N,q,\frac{L}{\nu},\frac{L_2}{\nu}$, $[u]_{C^{0,\lambda}(Q^+,\R^N)}$ and $\sigma$ and a positive number $\delta$ depending only on $n,N,\frac{L}{\nu},\frac{L_2}{\nu}$ and $[u]_{C^{0,\lambda}(Q^+,\R^N)}$ such that there holds
\begin{align*}
\lefteqn{ \Big(\mI{Z_{\sigma \rho}(x_0)} \kabs{Dv_h}^q \dx \Big)^{\frac{1}{q}}} \\
	& \leq \, c \, \Big[ \Big(\mI{Z_{(1-8 \sigma) \rho/2}(x_0)} 
	\! \kabs{Dv_h}^2 \dx \Big)^{\frac{1}{2}}
	+ \Big(\mI{Z_{(1-8 \sigma) \rho/2}(x_0)} \! 
	\big( 1 + \kabs{Du(x)}^2 + \kabs{Du(x + h e)}^2 \big)^{\frac{q}{2}} 
	\dx \Big)^{\frac{1}{q}}\Big] \\
	& \leq \, c \, \Big[ \kabs{h}^{-\frac{\alpha \lambda}{2}} 
	\Big(\mI{Z_{(1-8 \sigma) \rho/2}(x_0)} \kabs{\tau_{e,h}Du}^2 \dx \Big)^{\frac{1}{2}} 
	+ \Big(\mI{Z_{\rho/2}(x_0)} \big( 1 + \kabs{Du(x)}^2 \big)^{\frac{q}{2}} 
	\dx \Big)^{\frac{1}{q}}\Big] \\
	& \leq \, c \, \Big[
	\Big(\mI{Z_{\rho/2}(x_0)} \big( 1 + \kabs{Du}^2 \big) \dx \Big)^{\frac{1}{2}} 
	+ \Big(\mI{Z_{\rho/2}(x_0)} \big( 1 + \kabs{Du}^2 \big)^{\frac{q}{2}} 
	\dx \Big)^{\frac{1}{q}}\Big]
\end{align*}
for all $q \in [2,2+\delta)$. Here, we have also used the bound $\kabs{h}<2\sigma \rho$ (with $\kabs{\sigma} < \frac{1}{10}$) and the estimate \eqref{rbp-k-start} combined with the H\"older continuity of $u$ with exponent $\lambda$ (note that as a consequence the constant $c$ then depends additionally on the radius $\rho$). Hence, for all $s \in (2, \min\{s_0,2+\delta\})$ the previous inequality holds true. Keeping in mind the definition of $v_h$ and the higher integrability result \eqref{rbp-manni-higher-anfang}, we finally arrive at
\begin{equation*}
\Big(\mI{Z_{\sigma \rho}(x_0)} \kabs{\tau_{e,h}Du}^{s} \dx 
	\Big)^{\frac{1}{s}}	
	\, \leq \, c \, \, \kabs{h}^{\frac{\alpha \lambda}{2}} \, 
	\Big(\mI{Z_{\rho}(x_0)} \big( 1 + \kabs{Du}^2 
	\big) \dx \Big)^{\frac{1}{2}} \,,
\end{equation*}
which finishes the proof of the proposition.
\end{proof}
\end{proposition}

Moreover, we mention two direct consequences of Proposition~\ref{rbp-manni-diff-higher-anfang}. The first one follows from Lemma~\ref{rbp-manni-slice-lemma} and concerns the slice-wise mean-square deviation of $Du$:

\begin{corollary}
\label{rbp-manni-bem-slice}
Let $u \in W^{1,2}_{\Gamma}(Q_2^+,\R^N) \cap L^{\infty}(Q_2^+,\R^N) \cap C^{0,\lambda}(Q^+,\R^N)$ be a weak solution of \eqref{rbp-manni-model} under the assumptions (H1)-(H3) and (B). Furthermore, let  $Z_{\rho}(x_0) \subset Q^+$ for some $x_0 \in Q^+ \cup Q^0$ and $\sigma \in (0,\frac{1}{10})$. Then for every $\gamma \in (0,1)$ there exists a function $F_1 \in L^s(Z_{\sigma \rho}(x_0))$ such that the following estimate holds true:
\begin{multline*}
\qquad \Big( \mI{Z_{r}(z)} \babs{Du(x) - (Du)_{z',r}(x_n)}^2 \dx \Big)^{\frac{1}{2}} \\
\leq \, \Big( \mI{Z_{r}(z)} \mI{D_{r}(z')} \kabs{Du(x',x_n) - Du(y',x_n)}^2
	\dy' \dx \Big)^{\frac{1}{2}} 
\, \leq \, c \, r^{\frac{\gamma \alpha \lambda}{2}} \, F_1(z) \qquad
\end{multline*}
for all cylinders $Z_{r}(z) \subset Z_{\sigma \rho}(x_0)$ with $z \in Q^+ \cup Q^0$, and the constant $c$ depends only on $n,\alpha,\lambda$ and $\gamma$.
\end{corollary}

\begin{remarkon}
The $L^s$-norm of $F_1$ might blow up if $\gamma \nearrow 1$ (as a consequence of the application of the $L^q$-inequality for the maximal operator in the proof of Lemma~\ref{rbp-manni-slice-lemma}). Moreover, when verifying the assumptions of Lemma~\ref{rbp-manni-slice-lemma}, we observe that the number $K$ (resulting from the inequality in Proposition~\ref{rbp-manni-diff-higher-anfang}) depends on the radius $\rho$ and on $\sigma$. This dependency is reflected only in the $L^s$-norm of $F_1$. However, this will not be of importance because $\rho$ and $\sigma$ may be chosen fixed in every step of the subsequent iteration. More precisely, in the next section we will infer appropriate fractional Sobolev estimates on the cylinders $Z_{\sigma \rho}(x_0)$ and then, via a covering argument, also on $Q^+$ (respectively on smaller half-cubes in the course of the iteration).
\end{remarkon}

As a second consequence of Proposition~\ref{rbp-manni-diff-higher-anfang} we obtain a fractional Sobolev estimate for the tangential derivative $D'u$. This follows immediately from Lemma~\ref{frac-Sobol-lemma-2} and the inclusion $W^{\theta,s} \subseteq M^{\theta,s}$ (for $\theta \in (0,1)$, $s \in (1,\infty)$) given in Remarks~\ref{remark-frac-char}.

\begin{corollary}
\label{corollary-tangential-fractional}
Let $u \in W^{1,2}_{\Gamma}(Q_2^+,\R^N) \cap L^{\infty}(Q_2^+,\R^N) \cap C^{0,\lambda}(Q^+,\R^N)$ be a weak solution of \eqref{rbp-manni-model} under the assumptions (H1)-(H3) and (B). Then for every $\gamma \in (0,1)$ there holds 
\begin{equation*}
D'u = (D_1 u,\ldots,D_{n-1} u) \in M^{\gamma \alpha \lambda / 2,s}(Q_{\rho}^+,\R^{(n-1)N})
\end{equation*}
for every $\rho < 1$. In particular, there exists a function $H_1 \in L^s(Q_{1/2}^+)$ such that
\begin{equation*}
\kabs{D'u(x) - D'u(y)} \, \leq \, \kabs{x-y}^{\frac{\gamma \alpha \lambda}{2}} \, \big( H_1(x) + H_2(y) \big)
\end{equation*}
for almost all $x,y \in Q_{1/2}^+$.
\end{corollary} 

\subsection{A first estimate for the full derivative}
\label{first-full}

So far, we can estimate finite differences close to the boundary only with respect to tangential directions.  In order to find a fractional Sobolev estimate of type \eqref{frac-equiv-lemma} also with respect to the normal direction we next choose a cylinder $Z_{\rho}(x_0) \subset Q^+$, $x_0 \in Q^+ \cup Q^0$, $\rho \leq \widetilde{\rho}_{cacc}$ where $\widetilde{\rho}_{cacc}$ is from Lemma~\ref{manni-caccio-rev}, and $\sigma \in (0,\frac{1}{10})$. Furthermore, we fix a number $\gamma \in (0,1)$ to be specified later. We now study the model system \eqref{rbp-manni-model} on cylinders $Z_{r}(z)$ with $z \in Q^+ \cup Q^0$ such that $Z_{2r}(z) \subset Z_{\sigma \rho}(x_0)$, and by $M^*$ we always denote the maximal operator restricted to the cylinder $Z_{\sigma \rho}(x_0)$, i.\,e. 
\begin{equation*}
M^*(f)(z) \, := \, \sup_{Z_{\tilde{r}}(\tilde{z}) \subseteq Z_{\sigma \rho}(x_0), 
	\, z \in Z_{\tilde{r}}(\tilde{z})}
	\mI{Z_{\tilde{r}}(\tilde{z})} \kabs{f(x)} \dx \,.
\end{equation*}
for every $f \in L^1(Z_{\sigma \rho}(x_0),\R^k)$, $k \geq 1$, and $z \in Z_{\sigma \rho}(x_0)$. We shall frequently use the fact that the maximal operator is bounded as a mapping from $L^p$ to itself for every $p>1$.

\subsubsection*{A fractional Sobolev estimate for $\mathbf{a_n(\ccdot,u,Du)}$}

In coordinates we have the following representation of the weak formulation for the system \eqref{rbp-manni-model}:
\begin{equation*}
\sum_{j=1}^N \, \sum_{\kappa=1}^n \, \mI{Z_{r}(z)} a_{\kappa}^j (x,u(x),Du(x)) \, D_{\kappa} \varphi^{j} \dx
	\, = \, \sum_{j=1}^N \, \mI{Z_{r}(z)} b^j (x,u(x),Du(x)) \, \varphi^{j} \dx
\end{equation*}
for all $\varphi \in C_{0}^{\infty}(Z_{r}(z),\R^N)$. Following the approach of \cite{KRONZHABIL}, we are going to derive in the first step a weak differentiability result for the function 
\begin{equation}
\label{rbp-manni-def-A-r}
A_r^j(x_n) \, := \, \mI{D_r(z')} a_n^j(x',x_n,u(x',x_n),Du(x',x_n)) dx'
\end{equation}
for every $j \in \{1,\ldots,N\}$ and $x_n \in I_r(z_n)$. For this purpose we choose a ``splitting'' test function of the form $\varphi(x) = \phi_1(x') \, \phi_2(x_n) \, E_j$ where $\phi_1 \in C_0^{\infty}(D_r(z'))$ with $\phi_1 \equiv 1$ on the $(n-1)$-dimensional ball $D_{\tau r}(z')$ for some $\tau \in (0,1)$, $\phi_2 \in C_0^{\infty}(I_r(z_n))$, and where $E_j$ denotes the standard unit coordinate vector in $\R^N$. Employing the above identity with such a test function $\varphi$ then yields
\begin{align*}
\lefteqn{\hspace{-0.25cm} \mI{I_r(z_n)} \mI{D_r(z')} a_n^j(x,u(x),Du(x)) \, \phi_1(x') 
\, D_n \phi_2(x_n) \dx' \dx_n } \\
	& = \, - \mI{I_r(z_n)} \frac{1}{\kabs{D_r(z')}} \int_{D_r(z') \setminus D_{\tau r}(z')}
	\sum_{\kappa=1}^{n-1} \, a_{\kappa}^j (x,u(x),Du(x)) 
	\, D_{\kappa} \phi_1(x') \, \phi_2(x_n) \dx' \dx_n \\
	& \qquad {}+ \mI{I_r(z_n)} \mI{D_r(z')}  b^j (x,u(x),Du(x)) \, 
	\phi_1(x') \, \phi_2(x_n) \dx' \dx_n \\
	& = \, - \mI{I_r(z_n)} \frac{1}{\kabs{D_r(z')}} \int_{\tau r}^r 
	\int_{\partial D_{\widetilde{r}}(z')}
	\sum_{\kappa=1}^{n-1} \, \big[ a_{\kappa}^j (x,u(x),Du(x)) - 
	a_{\kappa}^j (z,(u)_{z,r},(Du)_{z',r}(x_n)) \big] \\
	& \hspace{5cm} \times D_{\kappa} \phi_1(x') \, d\Hm^{n-2}(x') 
	\, d\widetilde{r} \, \phi_2(x_n) \dx_n \\
	& \qquad {}+ \mI{I_r(z_n)} \mI{D_r(z')}  b^j (x,u(x),Du(x)) \, 
	\phi_1(x') \dx'  \, \phi_2(x_n) \dx_n
\end{align*}
for $j \in \{1,\ldots,N\}$, where we have used the co area formula in the last line. In particular, we may choose by approximation a cut-off function of the form
\begin{equation*}
\phi_1(x') \, = \,
\begin{cases} \begin{array}{c l} 
1 & \qquad \text{if } \kabs{x'-z'} \leq \tau r \,, \\
\frac{r-\kabs{x'-z'}}{(1-\tau) r} & \qquad \text{if } \tau r < \kabs{x'-z'} < r \,,  \\
0 & \qquad \text{if } \kabs{x'-z'} \geq r \,.
\end{array}
\end{cases}
\end{equation*}
We note that this implies $D_{\kappa} \phi_1(x') = - \frac{1}{(1-\tau)r} \, \frac{x_{\kappa} - z_{\kappa}}{\kabs{x'-z'}}$ for every $\kappa \in \{1,\ldots,n-1\}$ provided that $\kabs{x'-z'} \in (\tau r,r)$. Setting
\begin{equation}
\label{def-B}
B^j_{\kappa}(x) \, = \, a^j_{\kappa} (x,u(x),Du(x)) - a^j_{\kappa} (z,(u)_{z,r},(Du)_{z',r}(x_n)) 
\end{equation}
for $j \in \{1,\ldots,N\}$ and $\kappa \in \{1,\ldots,n-1\}$, we calculate with this particular choice for $\phi_1$:
\begin{align*}
\lefteqn{\hspace{-0.75cm} \mI{I_r(z_n)} \mI{D_r(z')} a_n^j(x,u(x),Du(x)) \, \phi_1(x') 
\dx' \, D_n \phi_2(x_n)  \dx_n } \\
& = \, \mI{I_r(z_n)} \frac{1}{\kabs{D_r(z')}} \mI{\tau r}^r \int_{\partial D_{\widetilde{r}}(z')}
	B^j(x) \cdot 
	\frac{x' - z'}{\kabs{x'-z'}} \, d\Hm^{n-2}(x') \, d\widetilde{r} \, \phi_2(x_n) \dx_n \\
	& \qquad {}+ \mI{I_r(z_n)} \mI{D_r(z')}  b^j (x,u(x),Du(x)) \, 
	\phi_1(x')  \dx' \, \phi_2(x_n) \dx_n \,.
\end{align*}
Recalling the definition of $A_r^j(x_n)$ given in \eqref{rbp-manni-def-A-r}, we consider the limit $\tau \nearrow 1$ and conclude from Lebesgue's differentiation Theorem that for almost every radius $r$ (and fixed center $z \in Z_{\sigma \rho}(x_0)$) such that $Z_{r}(z) \subset Z_{\sigma \rho}(x_0)$ there holds
\begin{align*}
\int_{I_r(z_n)} \!\! A_r^j(x_n) \, D_n \phi_2(x_n)  \dx_n
	& = \, \int_{I_r(z_n)} \! \frac{1}{\kabs{D_r(z')}} \int_{\partial D_{r}(z')} \!
	B^j(x) \cdot 
	\frac{x' - z'}{\kabs{x'-z'}} \, d\Hm^{n-2}(x') \, \phi_2(x_n) \dx_n \\
	& \qquad {}+ \int_{I_r(z_n)} \mI{D_r(z')}  b^j (x,u(x),Du(x))
	\dx' \, \phi_2(x_n) \dx_n \,.
\end{align*}
Hence, for almost every radius $r$ with $Z_{r}(z) \subset Z_{\sigma \rho}(x_0)$ we find that $A_r(x_n)=(A_r^1(x_n),\ldots,A_r^N(x_n))$ is weakly differentiable on $I_r(z_n)$ (note that the index $j \in \{1,\ldots,N\}$ and the test function $\phi_2$ are arbitrary in the latter identity), and its weak derivative is given by
\begin{equation}
\label{rbp-manni-A-abl}
A_r^{\prime}(x_n) \, = \, - \frac{1}{\kabs{D_r(z')}} \int_{\partial D_{r}(z')}
	B(x) \cdot \frac{x' - z'}{\kabs{x'-z'}} \, d\Hm^{n-2}(x') - \mI{D_r(z')}  b (x,u(x),Du(x)) \dx'\,.
\end{equation}
We next consider for any fixed $r$ all radii $\widetilde{\rho} \in (0,r]$ and we define the set $J$ via
\begin{equation*}
\label{def-J}
J \, = \, \Big\{ \widetilde{\rho} : \widetilde{\rho} \in (0,r] \quad \text{and} \quad
	\int_{I_{\widetilde{\rho}}(z_n)} \int_{\partial D_{\widetilde{\rho}}(z')} 
	\kabs{B(x)} \, d\Hm^{n-2}(x') \dx_n \, > \, \frac{2}{r} \int_{Z_r(z)} \kabs{B(x)} \dx \Big\} \,.
\end{equation*}
The following computations reveal that there holds $\Lm^1(J) < \frac{r}{2}$: employing the co area formula and Fubini's Theorem we get
\begin{align*}
\int_{Z_r(z)} \kabs{B(x)} \dx 
	& \geq \, \int_0^r \int_{I_{\widetilde{\rho}}(z_n)} \int_{\partial D_{\widetilde{\rho}}(z')} 
	\kabs{B(x)} \, d\Hm^{n-2}(x') \dx_n \, d\widetilde{\rho} \\
	& \geq \, \int_J \int_{I_{\widetilde{\rho}}(z_n)} \int_{\partial D_{\widetilde{\rho}}(z')} 
	\kabs{B(x)} \, d\Hm^{n-2}(x') \dx_n \, d\widetilde{\rho} \\
	& > \, \int_J \frac{2}{r} \int_{Z_r(z)} \kabs{B(x)} \dx \, d\widetilde{\rho} 
	\, = \, \Lm^1(J) \, \frac{2}{r} \int_{Z_r(z)} \kabs{B(x)} \dx \,.
\end{align*}
Therefore, we find some radius $\bar{\rho} \in [\frac{r}{2},r]$ such that on the one hand $A_{\bar{\rho}}(x_n)$ is weakly differentiable and on the other hand $\bar{\rho} \notin J$. 
Hence, in view of Poincar\'{e}'s inequality and identity \eqref{rbp-manni-A-abl}, we obtain for this choice of $\bar{\rho}$:
\begin{align}
\label{rbp-manni-A-MW}
\lefteqn{\hspace{-0.75cm} \mI{I_{\bar{\rho}}(z_n)} \babs{A_{\bar{\rho}}(x_n) - (A_{\bar{\rho}})_{z_n,\bar{\rho}}} \dx_n
	\, \leq \, c(N) \int_{I_{\bar{\rho}}(z_n)} \babs{A^{\prime}_{\bar{\rho}}(x_n)} \dx_n} \nonumber \\
	& \leq \, \frac{c(N)}{\kabs{D_{\bar{\rho}}(z')}} \int_{I_{\bar{\rho}}(z_n)} 
	\int_{\partial D_{\bar{\rho}}(z')} \kabs{B(x)} \, d\Hm^{n-2}(x') \dx_n \nonumber \\
	& \qquad {}+ c(N) \int_{I_{\bar{\rho}}(z_n)} \mI{D_{\bar{\rho}}(z')}  
	\kabs{b(x,u(x),Du(x))} \dx' \dx_n \nonumber \\
	& \leq \, c(N) \Big[ \frac{1}{\kabs{D_{\bar{\rho}}(z')} \, r} \int_{Z_r(z)} \kabs{B(x)} \dx 
	+ \bar{\rho} \mI{Z_{\bar{\rho}}(z)} \kabs{b(x,u(x),Du(x))} \dx \Big] \nonumber \\
	& \leq \, c(n,N) \Big[ \mI{Z_r(z)} \kabs{B(x)} \dx 
	+ r \mI{Z_r(z)} \kabs{b(x,u(x),Du(x))} \dx \Big] \,.
\end{align}
In the next step we control the integrals arising on the right-hand side of the last inequality by using the growth conditions on coefficients and inhomogeneity, respectively, and by exploiting the assumption that $u$ is H\"older continuous with exponent $\lambda$ (which shall be used without any further comment).

For the first integral in \eqref{rbp-manni-A-MW} we use the definition of $B(x)$ in \eqref{def-B}, the assumptions (H1), (H3), and Corollary~\ref{rbp-manni-bem-slice} to see
\begin{align*}
\mI{Z_r(z)} \babs{B(x)} \dx & \leq \, \mI{Z_r(z)} \big[ \, \babs{a(x,u(x),Du(x)) - 
	a(z,(u)_{z,r},Du(x))} \\
	& \hspace{2cm} {}+ \babs{a(z,(u)_{z,r},Du(x)) - 
	a(z,(u)_{z,r},(Du)_{z',r}(x_n))} \, \big] \dx \\
	& \leq \, 4 \, L \big( r^{\alpha} + [u]^{\alpha}_{C^{0,\lambda}(Q^+,\R^N)} 
	\, r^{\alpha \lambda}\big) \mI{Z_r(z)} \big( 1 + \kabs{Du} \big) \dx \\
	& \qquad {}+ L \mI{Z_r(z)} \babs{Du(x) - (Du)_{z',r}(x_n)} \dx  \\
	& \leq \, c \, r^{\frac{\gamma \alpha \lambda}{2}} \, 
	\big( M^*\big(1+\kabs{Du}\big)(z) +  F_1(z) \big) \,,
\end{align*}
and the constant $c$ depends only on $n,L,[u]_{C^{0,\lambda}(Q^+,\R^N)},\alpha,\lambda$ and $\gamma$. Moreover, the functions $F_1$ and $M^*\big(1+\kabs{Du}\big)$ belong to the space $L^s(Z_{\sigma \rho}(x_0))$, due to Corollary~\ref{rbp-manni-bem-slice} and the higher integrability of $Du$ (combined with standard properties of the maximal function).

For the second integral in \eqref{rbp-manni-A-MW}, we initially assume that we are close to the boundary, meaning that $z_n < 2r$. We then infer the following estimate from the natural growth condition (B), the Caccioppoli inequality from Lemma~\ref{manni-caccio-rev} (note that $2r \leq \widetilde{\rho}_{cacc}$), and the Poincar\'{e} inequality in the boundary version: 
\begin{align}
\label{rbp-manni-A-inhomogenitaet}
r \mI{Z_r(z)} \kabs{b(x,u(x),Du(x))} \dx & \leq \, r \, \mI{Z_r(z)} (L + L_2 \, \kabs{Du}^2) \dx \nonumber \\
	& \leq \, c \, \Big( r^{1-1 + \lambda} \mI{Z_{2r}(z)} \kabs{Du} \dx 
	+ r^{2\alpha +1}\Big) + r \, L \nonumber \\
	& \leq \, c \, r^{\lambda} \, M^*\big(1+\kabs{Du}\big)(z) \,,
\end{align}
and the constant $c$ depends only on $n,N,L,L_2,\nu$ and $[u]_{C^{0,\lambda}(Q^+,\R^N)}$. For cylinders in the interior, where $z_n \geq 2r$, we end up with exactly the same estimate using interior versions of Caccioppoli and Poincar\'{e} where $\kabs{u}$ is replaced by $\kabs{u - (u)_{z,2r}}$.

Hence, combining the last two estimates, we conclude from \eqref{rbp-manni-A-MW}
\begin{multline}
\label{rbp-manni-A-MW-2}
\quad \mI{Z_{\bar{\rho}}(z_n)} \Babs{\mI{D_{\bar{\rho}}(z')} a_n(y',x_n,u(y',x_n),Du(y',x_n)) \, dy' - 
	\mI{Z_{\bar{\rho}}(z)} a_n(\tilde{y},u(\tilde{y}),Du(\tilde{y})) \, d\tilde{y}} \dx \\
	= \, \mI{I_{\bar{\rho}}(z_n)} \babs{A_{\bar{\rho}}(x_n) - 
	(A_{\bar{\rho}})_{z_n,\bar{\rho}}} \dx_n 
	\, \leq \, c \, r^{\frac{\gamma \alpha \lambda}{2}} \, 
	\big[ M^*\big(1+\kabs{Du}\big)(z) + F_1(z)\big] \,, \quad
\end{multline}
and the constant $c$ depends only on $n,N,L,L_2,\nu,[u]_{C^{0,\lambda}(Q^+,\R^N)},\alpha,\lambda$ and $\gamma$. Besides, we have $F_1$, $M^*\big(1+\kabs{Du}\big) \in L^s(Z_{\sigma \rho}(x_0))$ for some $s > 2$. Furthermore, applying Jensen's inequality, conditions (H1), (H3), and Corollary~\ref{rbp-manni-bem-slice} we find
\begin{align}
\label{rbp-manni-A-MW-3}
\lefteqn{ \mI{Z_{\bar{\rho}}(z)} \Babs{a_n(x,u(x),Du(x)) - \mI{D_{\bar{\rho}}(z')} a_n(y',x_n,u(y',x_n),Du(y',x_n)) \, dy'} \dx } \nonumber \\
	& \leq \, c\big(L,[u]_{C^{0,\lambda}(Q^+,\R^N)}\big) \, \bar{\rho}^{\, \alpha \lambda}
	\mI{Z_{\bar{\rho}}(z)} \big(1+\kabs{Du}\big) \dx \nonumber \\
	& \qquad {} +  L \mI{Z_{\bar{\rho}}(z)} \mI{D_{\bar{\rho}}(z')} 
	\babs{Du(x',x_n) - Du(y',x_n)} \dy' \dx \nonumber \\
	& \leq \, c\big(n,L, [u]_{C^{0,\lambda}(Q^+,\R^N)},\alpha,\lambda,\gamma\big) 
	\, \bar{\rho}^{\,\frac{\gamma \alpha \lambda}{2}} \, 
	\big[ M^*\big(1+\kabs{Du}\big)(z) + F_1(z)\big] \,.
\end{align}
Combining \eqref{rbp-manni-A-MW-2} and \eqref{rbp-manni-A-MW-3}, we conclude
\begin{equation*}
\mI{Z_{\bar{\rho}}(z)} \babs{a_n(x,u(x),Du(x)) - \big( a_n(\ccdot,u,Du)\big)_{z,\bar{\rho}} } \dx \, \leq \, c \, r^{\frac{\gamma \alpha \lambda}{2}} \, 
	\big[ M^*\big(1+\kabs{Du}\big)(z) + F_1(z)\big] 
\end{equation*}
for every $r$ with $Z_{r}(z) \subset Z_{\sigma \rho}(x_0)$ and an appropriate radius $\bar{\rho} \in [\frac{r}{2},r]$ for which $A_{\bar{\rho}}(x_n)$ is weakly differentiable on $I_r(z_n)$ and $\bar{\rho} \notin J$. The constant $c$ here depends only on $n,N,L,L_2,\nu$, $[u]_{C^{0,\lambda}(Q^+,\R^N)},\alpha,\lambda$ and $\gamma$. In particular, this yields
\begin{equation*}
\mI{Z_{r/2}(z)} \babs{a_n(x,u(x),Du(x)) - \big( a_n(\ccdot,u,Du)\big)_{z,r/2} } \dx \, \leq \, c \, r^{\frac{\gamma \alpha \lambda}{2}} \, 
	\big[ M^*\big(1+\kabs{Du}\big)(z) + F_1(z)\big] \,,
\end{equation*}
and the constant $c$ admits the same dependencies as in the previous inequality. This allows to apply the characterization of fractional Sobolev spaces given in Lemma~\ref{lemma-frac-char} and Remarks~\ref{remark-frac-char} (note that these results also hold true if we replace the balls by cubes or cylinders). Since the cylinders $Z_{\rho}(x_0) \subset Q^+$ were chosen arbitrarily we infer via a covering argument
\begin{equation*}
a_n(\ccdot,u,Du) \in M^{\frac{\gamma \alpha \lambda}{2},s}(Q_{1/2}^+,\R^N) \,.
\end{equation*}
Furthermore, there exists a function $G_1 \in L^{s}(Q_{1/2}^+,\R^N)$ which satisfies
\begin{equation*}
\kabs{a_n(x,u(x),Du(x)) - a_n(y,u(y),Du(y))} \, \leq \, 
\kabs{x-y}^{\frac{\gamma \alpha \lambda}{2}} \, \big( G_1(x) + G_1(y) \big) 
\end{equation*}
for almost every $x,y \in Q_{1/2}^+$. We finally note that $G_1$ can be calculated from  $c$, $M^*\big(1+\kabs{Du}\big)$, $F_1(z)$ and the restriction on the radius $\rho$.

We close this first step with some remarks concerning the components $a_k(\cdot,u,Du)$ of the coefficients, $k \in \{1,\ldots,n-1\}$, and the interior situation:

\begin{remarks}
\label{bem-manni-innen}
We first note that testing the system \eqref{rbp-manni-model} with finite differences in normal direction of the weak solution $u$ is not allowed. Hence, the statement in Proposition~\ref{rbp-manni-diff-higher-anfang} cannot be expected to cover (via a modified proof) also differences of $Du$ in \emph{any} arbitrary direction $e \in S^{n-1}$ up to the boundary. This reveals the crucial point for the up-to-the-boundary estimates derived in this section: the method makes only an up to the boundary estimate for $a_n(\cdot,u,Du)$ available -- which is still sufficient to enable us later to find an appropriate fractional Sobolev estimate for $Du$~-- but a corresponding estimate for $a_k(\cdot,u,Du)$, $k \in \{1,\ldots,n-1\}$, does not follow.

For cylinders in the interior, however, Proposition~\ref{rbp-manni-diff-higher-anfang} holds true for every direction $e \in S^{n-1}$. As a consequence, we may repeat the arguments above line-by-line and end up with an interior fractional estimate for the full coefficients $a(\cdot,u,Du)$. We here mention that fractional Sobolev estimates for the coefficients $a(\cdot,u,Du)$ are not necessary in the interior to prove the dimension reduction for the singular set. In fact, interior fractional Sobolev estimates for weak solutions to elliptic systems with inhomogeneities obeying a natural growth condition can be obtained directly by exploiting the fundamental estimate \eqref{rbp-k-start}, see \cite{MINGIONE03}.
\end{remarks}

\subsubsection*{A fractional Sobolev estimate for $\mathbf{Du}$}

The ellipticity condition (H2) and the upper bound in (H1) allow to estimate
\begin{align*}
\lefteqn{\hspace{-0.25cm} \big[ a_n(x,u(x),Du(x)) - a_n(x,u(x),Du(y))\big] \cdot \big( D_n u(x) - D_n u(y) \big)} \\
	& = \, \int_0^1 D_z a_n\big(x,u(x),Du(y) + t(Du(x)-Du(y))\big) \dt \\[-0.2cm]
	& \hspace{3cm} \big( D u(x) - D u(y) \big) \cdot \big( D_n u(x) - D_n u(y) \big) \\
	& \geq \, \nu \, \kabs{D_n u(x) - D_n u(y)}^2 
	- L \, \kabs{D' u(x) - D' u(y)} \, \kabs{D_n u(x) - D_n u(y)}
\end{align*}
for almost all $x,y \in Q_{1/2}^+$. Dividing by $\kabs{D_n u(x) - D_n u(y)}$ (provided that $D_n u(x) \neq D_n u(y)$ which is the nontrivial case) and taking into account the fractional Sobolev estimates for both $a_n(\cdot,u,Du)$ and the tangential derivative $D'u$ from Corollary~\ref{corollary-tangential-fractional} and condition (H3), the latter inequality implies
\begin{align*}
\label{D_n_u_inherit}
\lefteqn{\hspace{-0.75cm} \nu \,  \kabs{D_n u(x) - D_n u(y)} 
	\, \leq \, \babs{a_n(x,u(x),Du(x)) - a_n(x,u(x),Du(y))} + L \, \kabs{D' u(x) - D' u(y)} } \\
	& \leq \, L \, \big( \kabs{x-y}^{\alpha} + [u]^{\alpha}_{C^{0,\lambda}(Q^+,\R^N)} 
	\kabs{x-y}^{\alpha \lambda} \big) \, \big( 1 + \kabs{Du(y)} \big) \\
	& \qquad {} + \kabs{x-y}^{\frac{\gamma \alpha \lambda}{2}} \, 
	\big( G_1(x) + G_1(y) \big) + L \, \kabs{x-y}^{\frac{\gamma \alpha \lambda}{2}} \, 
	\big( H_1(x) + H_1(y) \big)  \\
	& \leq \, c(L,[u]_{C^{0,\lambda}(Q^+,\R^N)}) \, \kabs{x-y}^{\frac{\gamma \alpha \lambda}{2}}
	\, \big( 1 + \kabs{Du(y)} + G_1(x) + G_1(y) + H_1(x) + H_1(y)  \big) 
\end{align*}
for almost every $x,y \in Q_{1/2}^+$, meaning that we have $D_n u \in M^{\frac{\gamma \alpha \lambda}{2},s}(Q_{1/2}^+,\R^N)$. Combined with Corollary~\ref{corollary-tangential-fractional} we hence end up with 
\begin{equation*}
D u \in M^{\frac{\gamma \alpha \lambda}{2},s}(Q_{1/2}^+,\R^{nN}) \,,
\end{equation*}
which is the desired estimate for the full derivative $Du$. We recall the embedding for the fractional Sobolev spaces, namely that
\begin{equation*}
M^{\gamma \alpha \lambda/2,s}(Q_{1/2}^+,\R^{nN}) \subset W^{\gamma' \gamma \alpha \lambda/2,s}(Q_{1/2}^+,\R^{nN})
\end{equation*}
for all $\gamma' \in (0,1)$. Then, in view of the interpolation Theorem~\ref{CampInterpolation} and the fact that $\gamma$ and $\gamma'$ may be chosen arbitrarily close to $1$ (an appropriate choice is for example $\gamma=\gamma'=(\tfrac{n}{n+2\lambda})^{1/2}$), we finally arrive at the higher integrability result 
\begin{equation*}
Du \in L^{s(1+\alpha \lambda / 2)}(Q_{1/2}^+,\R^{nN}) \, .
\end{equation*}

\subsection{Iteration}
\label{sec-iteration}

In the next step we iterate the fractional Sobolev estimate for $Du$. To this aim we define a sequence $(b_k)_{k \in \N}$ as follows:
\begin{equation*}
b_0 \, := \, 0, \qquad b_{k+1} \, := \, \frac{\alpha \lambda}{2} + b_k \, \Big(1-\frac{\lambda}{2}\Big) \, = \, b_k + \frac{\lambda}{2} \, ( \alpha - b_k )
\end{equation*}
for all $k \in \N_0$. We observe that the sequence $(b_k)$ is increasing with $b_k \nearrow \alpha$. The strategy of the proof is the following: For every $k \in \N_0$ we show by induction the following inclusions:
\begin{align*}
Du \in L^{s_{k}(1+b_{k})}(Q_{2^{-k}}^+,\R^{nN}) \quad & \rightarrow \quad Du \in M^{\gamma b_{k+1},s_{k+1}}(Q_{1/2^{k+1}}^+,\R^{nN}) \\
	& \rightarrow \quad Du \in L^{s_{k+1}(1+b_{k+1})}(Q_{1/2^{k+1}}^+,\R^{nN}) \,,
\end{align*}
where $\gamma \in (0,1)$ is an arbitrary number and where $(s_k)_{k \in \N}$ is a decreasing sequence of higher integrability exponents with $s_k > 2$ for every $k \in \N_0$. The first step of the induction, $k=0$, was already performed above (with $s_1=s$). We now proceed to the inductive step: The objective is to find the first inclusion by improving the fractional Sobolev estimates in Sections~\ref{high-anfang} and~\ref{first-full}, and then to deduce in the second step the higher integrability result by applying the interpolation Theorem~\ref{CampInterpolation}.

\subsubsection*{Higher integrability II}

We again need to deduce a higher integrability result for the tangential differences $\tau_{e,h}Du$ (cf. Proposition~\ref{rbp-manni-diff-higher-anfang}) which now incorporates the fact that $Du$ is assumed to be integrable with exponent $s_k(1 + b_k)$. In what follows we will frequently use a simple consequence of $b_k \leq \alpha$, namely the inequality
\begin{equation*}
\alpha \lambda + b_k (1-\lambda) \, \geq \, \frac{\alpha \lambda}{2} + b_k \, \Big(1-\frac{\lambda}{2}\Big) \, = \,  b_{k+1} \,.
\end{equation*}

\begin{proposition}
\label{rbp-manni-diff-higher-iteration}
Let $u \in  W^{1,2}_{\Gamma}(Q_2^+,\R^N) \cap L^{\infty}(Q_2^+,\R^N) \cap C^{0,\lambda}(Q^+,\R^N)$ be a weak solution to \eqref{rbp-manni-model} under the assumptions (H1)-(H3) and (B). Assume further $u \in W^{1,s_k (1 + b_k)}_{\Gamma}(Q_{2^{-k}}^+,\R^N)$ for some $k \in \N$, $s_k > 2$, and let $Z_{\rho}(x_0) \subset Q_{2^{-k}}^+$ for some $x_0 \in Q^0_{2^{-k}} \cup Q_{2^{-k}}^+$, $\sigma \in (0,\frac{1}{5})$, $e \in S^{n-1}$ with $e \perp e_n$ and $h \in \R$ satisfying $\kabs{h} \in (0,2\sigma\rho)$. Then there exists a higher integrability exponent $s_{k+1} \in (2,s_k)$ depending only on $n,N,\frac{L}{\nu},\frac{L_2}{\nu}$ and $[u]_{C^{0,\lambda}(Q^+,\R^N)}$ such that
\begin{equation*}
\mI{Z_{\sigma \rho}(x_0)} \kabs{\tau_{e,h}Du}^{s_{k+1}} \dx 
	\, \leq \, c \, \, \kabs{h}^{s_{k+1} b_{k+1}} \, 
	\Big(\mI{Z_{\rho}(x_0)} \big( 1 + \kabs{Du(x)} 
	\big)^{s_k(1+b_k)} \dx \Big)^{\frac{s_{k+1}}{s_k}}
\end{equation*}
for a constant $c= c\big(n,N,\frac{L}{\nu},\frac{L_2}{\nu},[u]_{C^{0,\lambda}(Q^+,\R^N)},\rho,\sigma\big)$.

\begin{proof}
We give only a sketch of proof and refer to \cite[Proof of Proposition 8.8]{BECK08} for more details. We start from the preliminary estimate \eqref{rbp-k-start} and show that for every $\theta \in (0,1)$ and every cylinder $Z_{r}(x_0) \subset Q_{2^{-k}}^+$ there holds
\begin{equation}
\label{rbp-manni-diff-iteration}
\int_{Z_{\theta r}(x_0)} \kabs{\tau_{e,h} Du}^2 \dx \, \leq \, c \, \kabs{h}^{2 b_{k+1}} \int_{Z_{r}(x_0)} \big( 1 + \kabs{Du} \big)^{2 + 2 b_k} \dx
\end{equation}
for all $e \in S^{n-1}$ with $e \perp e_n$, $h \in \R$ satisfying $\kabs{h} < \frac{r(1-\theta)}{2}$ and a constant $c$ depending only on $n,N,\tfrac{L}{\nu},\tfrac{L_2}{\nu} ,[u]_{C^{0,\lambda}(Q^+,\R^N)}, \theta$ and $r$. For this purpose, a suitable cut-off function is chosen, and the different terms arising on the right-hand side of \eqref{rbp-k-start} are then estimated taking advantage of standard properties of finite differences, the integrability of $Du$ with exponent $2 +  2 b_k$ and the H\"older continuity with exponent $\lambda$, see also \cite[p. 387]{MINGIONE03}.

In the next step we proceed similarly to the case $k=0$ and estimate the $L^{s_{k+1}}$-norm of $\kabs{\tau_{e,h} Du}$ for some exponent $s_{k+1}>2$ in terms of an appropriate power of $\kabs{h}$. To this end we consider directions $e \in S^{n-1}$ with $e \perp e_n$ and $h \in \R$ satisfying $\kabs{h}<2^{-k}$; furthermore, analogously to the proof of Proposition~\ref{rbp-manni-diff-higher-anfang} we set
\begin{equation*}
v_{h,k} \, := \, \frac{\tau_{e,h}u}{\kabs{h}^{b_{k+1}}} \,, \hspace{1cm}
\widetilde{{\cal A}}_{k}(h) \, := \, \frac{-{\cal A} (h)}{\kabs{h}^{b_{k+1}}} \,, \hspace{1cm}
\widetilde{{\cal B}}_{k}(h) \, := \, \frac{-{\cal B} (h)}{\kabs{h}^{b_{k+1}}} \,,
\end{equation*}
and $\widetilde{{\cal C}}_{k}(h) = \widetilde{{\cal C}}(h) = \int_0^1 D_z a\big(x,u(x),Du(x) + t\tau_{e,h} Du(x))\big) \, dt$ as above. Analogously to the derivation of \eqref{system-v-h} we then see that the map $v_{h,k} \in W^{1,2+2b_k}(Q^+_{2^{-k}-\kabs{h}},\R^N)$ is a weak solution to a linear system, for which the various terms need to be estimated in terms of the $L^2$-norms of $v_{h,k}$ and $Dv_{h,k}$. The only point differing from the estimates before is the one involving $\widetilde{{\cal B}}^{(k)}(h)$: to find an adequate inequality we first take advantage of the H\"older continuity of $u$ and Young's inequality and we see
\begin{align*}
\lefteqn{\hspace{-0.75cm} \int_{B^+_{R}(y)} \big( 1 + \kabs{Du(x + h e)} \big)^2 \, \kabs{\tau_{e,h}u}^{2 \alpha} \dx} \\
	& \leq \,  c\big([u]_{C^{0,\lambda}(B^+,\R^N)}\big) \, \kabs{h}^{2 \alpha \lambda - 
	2 b_k \lambda} \int_{B^+_{R}(y)} 
	\big( 1 + \kabs{Du(x + h e)} \big)^2 \, \kabs{\tau_{e,h}u}^{2 b_k} \dx \\
	& \leq \, c\big([u]_{C^{0,\lambda}(B^+,\R^N)}\big) \,\kabs{h}^{2 b_{k+1}} 
	\int_{B^+_{R}(y)} \big( 1 + \kabs{Du(x + h e)} + \kabs{G_h(x)} \big)^{2 + 2 b_k} \dx \,.
\end{align*}
Here we have used the fact that 
\begin{equation*}
\kabs{\tau_{e,h}u} \, \leq \, \kabs{h} \int_0^1 \kabs{Du(x+the)} \dt \, =: \kabs{h} \, G_h(x) \,,
\end{equation*}
and the function $G_h$ is $L^{s_k(1+b_k)}$-integrable on $B^+_R(y)$ in view of Fubini's Theorem:
\begin{equation*}
\int_{B^+_R(y)} \kabs{G_h}^{s_k (1+b_k)} \dx \, \leq \, \int_{Q^+_{2^{-k}}} \kabs{Du}^{s_k(1+ b_k)} \dx \, < \, \infty\,.
\end{equation*}
Hence, we find with Young's inequality for every $\epsilon \in (0,1)$
\begin{itemize}
 \item $\displaystyle \int_{B_R^+(y)} \babs{\widetilde{{\cal B}}^{(k)} (h) \cdot D\varphi} \dx
	\, \leq \, L \int_{B_R^+(y)} \kabs{h}^{-b_{k+1}} \,  \big( 1 + \kabs{Du(x + h e)} \big)
	\, \kabs{\tau_{e,h}u}^{\alpha} \, \kabs{D\varphi} \dx $ \\[0.1cm]
	$\displaystyle \hspace{1cm} \leq \, \epsilon \int_{B_R^+(y)} \eta_k^2 \, \kabs{Dv^{(k)}_h}^2 \dx 
	+  \frac{c \, \epsilon}{R^2} \int_{B_R^+(y)} \kabs{v^{(k)}_h}^2 \dx$ \\[0.1cm]
	$\displaystyle \hspace{1.75cm} {}+ c\big([u]_{C^{0,\lambda}(B^+,\R^N)}\big) \, 
	\epsilon^{-1} L^2 \, \int_{B_R^+(y)} 
	\big( 1 + \kabs{Du(x + h e)} + \kabs{G_h(x)}\big)^{2+2b_k} \dx$.
\end{itemize}
Arguing exactly as in the proof of Proposition~\ref{rbp-manni-diff-higher-anfang}, we obtain via the Sobolev-Poincar\'{e} inequality a reverse H\"older-type inequality, from which (taking advantage of the higher integrability of $G_h$) we then deduce the desired assertion by the global Gehring Lemma.
\end{proof}
\end{proposition}

\begin{remark}
\label{rem-limit}
If $s_k (1 + b_k) \geq 2 + 2 \alpha$ is satisfied, it is easy to check that the inequality \eqref{rbp-manni-diff-iteration} and in turn the statement of the proposition on the higher integrability of the differences hold true with $b_{k+1}$ replaced by $\alpha$.
\end{remark}

Proposition~\ref{rbp-manni-diff-higher-iteration} combined with Lemma~\ref{rbp-manni-slice-lemma} and with Lemma~\ref{frac-Sobol-lemma-2}, respectively, again allows to state two direct consequences concerning the slice-wise mean-square deviation of $Du$ and a suitable fractional differentiability of the tangential derivative $D'u$:

\begin{corollary}
\label{rbp-manni-bem-slice-iteration}
Let $u \in  W^{1,2}_{\Gamma}(Q_2^+,\R^N) \cap L^{\infty}(Q_2^+,\R^N)\cap C^{0,\lambda}(Q^+,\R^N)$ be a weak solution to \eqref{rbp-manni-model} under the assumptions (H1)-(H3) and (B). Assume further $u \in W^{1,s_k (1 + b_k)}_{\Gamma}(Q_{2^{-k}}^+,\R^N)$ for some $k \in \N$, $s_k > 2$, and let $Z_{\rho}(x_0) \subset Q_{2^{-k}}^+$ for some $x_0 \in Q^0_{2^{-k}} \cup Q_{2^{-k}}^+$ and $\sigma \in (0,\frac{1}{5})$. Then for every $\gamma \in (0,1)$ there exists a function $F_{k+1} \in L^{s_{k+1}}(Z_{\sigma \rho}(x_0))$ where $s_{k+1}\in (2,s_k)$ is the higher integrability exponent determined in Proposition~\ref{rbp-manni-diff-higher-iteration} such that the following estimate holds true:
\begin{multline*}
\qquad \Big( \mI{Z_{r}(z)} \babs{Du(x) - (Du)_{z',r}(x_n)}^2 \dx \Big)^{\frac{1}{2}} \\
\leq \, \Big( \mI{Z_{r}(z)} \mI{D_{r}(z')} \kabs{Du(x',x_n) - Du(y',x_n)}^2
	\dy' \dx \Big)^{\frac{1}{2}} 
\, \leq \, c \, r^{\gamma b_{k+1}} \, F_{k+1}(z) \qquad
\end{multline*}
for all cylinders $Z_{r}(z) \subset Z_{\sigma \rho}(x_0)$ with $z \in Q^+ \cup Q^0$, and the constant $c$ depends only on $n,\alpha,\lambda$ and $\gamma$.
\end{corollary}

\begin{corollary}
\label{corollary-tangential-fractional-iteration}
Let $u \in W^{1,2}_{\Gamma}(Q_2^+,\R^N) \cap L^{\infty}(Q_2^+,\R^N) \cap C^{0,\lambda}(Q^+,\R^N)$ be a weak solution to \eqref{rbp-manni-model} under the assumptions (H1)-(H3) and (B). Assume further $u \in W^{1,s_k (1 + b_k)}_{\Gamma}(Q_{2^{-k}}^+,\R^N) $ for some $k \in \N$, $s_k > 2$. Then for every $\gamma \in (0,1)$ there holds 
\begin{equation*}
D'u \in M^{\gamma b_{k+1},s_{k+1}}(Q_{\rho}^+,\R^{(n-1)N})
\end{equation*}
for every $\rho < \frac{1}{2^{k+1}}$. In particular, there exists a function $H_{k+1} \in L^{s_{k+1}}(Q_{1/2^{k+1}}^+)$ such that
\begin{equation*}
\kabs{D'u(x) - D'u(y)} \, \leq \, \kabs{x-y}^{\gamma b_{k+1}} \, \big( H_{k+1}(x) + H_{k+1}(y) \big)
\end{equation*}
for almost all $x,y \in Q_{1/2^{k+1}}^+$.
\end{corollary} 

\subsubsection*{An improved fractional Sobolev estimate for $\mathbf{a_n(\ccdot,u,Du)}$}

Taking into account that $Du$ is assumed to be integrable with exponent $s_k(1 + b_k)$, we next proceed similarly to the case $k=0$: We choose a cylinder $Z_{\rho}(x_0) \subset Q^+_{2^{-k}}$ with center $x_0 \in Q^+_{2^{-k}} \cup Q^0_{2^{-k}}$ and radius $\rho$ sufficiently small , i.\,e. $\rho \leq \widetilde{\rho}_{cacc}$ where $\widetilde{\rho}_{cacc}$ is from the Caccioppoli-type inequality in Lemma~\ref{manni-caccio-rev}, and $\sigma \in (0,\frac{1}{5})$. Furthermore, we fix a number $\gamma \in (0,1)$ and again study the model system \eqref{rbp-manni-model} on cylinders $Z_{r}(z)$ with $z \in Q^+_{2^{-k}} \cup Q^0_{2^{-k}}$ such that $Z_{2r}(z) \subset Z_{\sigma \rho}(x_0)$. Using the notation from Section~\ref{first-full}, we first improve the estimate \eqref{rbp-manni-A-MW-2}: To this aim we start with inequality \eqref{rbp-manni-A-MW}: For the first integral on the right-hand side of \eqref{rbp-manni-A-MW} we recall the definition of $B(x)$ in \eqref{def-B} and take advantage of conditions (H1) and (H3) to infer
\begin{align*}
\mI{Z_r(z)} \babs{B(x)} \dx 
	& \leq \, L \mI{Z_r(z)} \big( \kabs{x-z}^{\alpha} + \kabs{u(x)-(u)_{z,r}}^{\alpha} 
	\big) \, \big( 1 + \kabs{Du(x)} \big) \dx  \\
	& \qquad {}+ L \mI{Z_r(z)} \babs{Du(x) - (Du)_{z',r}(x_n)} \dx \,.
\end{align*}
In view of H\"older's and Jensen's inequality, the H\"older continuity of $u$ and Poincar\'e's Lemma, we derive
\begin{align}
\label{u-Du-iteration}
 \lefteqn{\hspace{-0.75cm} \mI{Z_r(z)}  \kabs{u(x)-(u)_{z,r}}^{\alpha} \,
	\big( 1 + \kabs{Du(x)} \big) 
	\dx} \nonumber \\
	& \leq \, \Big( \mI{Z_r(z)}  \kabs{u(x)-(u)_{z,r}}^{\alpha \frac{1+b_k}{b_k}} \dx 
	\Big)^{\frac{b_k}{1 + b_k}}  \Big( \mI{Z_r(z)} \big( 1 + \kabs{Du} )^{1+b_k} \dx 
	\Big)^{\frac{1}{b_k + 1}} \nonumber \\
	& \leq \, c \, r^{\alpha \lambda - b_k \lambda} 
	\,  \Big( \mI{Z_r(z)} \kabs{u(x)-(u)_{z,r}}^{1 + b_k} \dx 
	\Big)^{\frac{b_k}{1+b_k}}  \Big( \mI{Z_r(z)} \big( 1 + \kabs{Du} )^{1+b_k} \dx 
	\Big)^{\frac{1}{1+b_k}} \nonumber \\
	& \leq \, c \, r^{\alpha \lambda + b_k (1- \lambda)} 
	\mI{Z_r(z)} \big( 1 + \kabs{Du} )^{1+b_k} \dx 
	\, \leq \, c \, r^{\gamma b_{k+1}} \,  M^*\big((1+\kabs{Du})^{1+b_k}\big)(z)
\end{align}
for $c=c(n,[u]_{C^{0,\lambda}(Q^+,\R^N)})$. Furthermore, we trivially have 
\begin{align*}
\mI{Z_r(z)} \kabs{x-z}^{\alpha} \, \big( 1 + \kabs{Du(x)} \big) \dx
	& \leq \, c(n) \, r^{\gamma b_{k+1}} \,  M^*\big((1+\kabs{Du})^{1+b_k}\big)(z) \,.
\end{align*}
Keeping in mind Corollary~\ref{rbp-manni-bem-slice-iteration} we finally arrive at the following estimate for the integral of $\kabs{B(x)}$:
\begin{equation}
\label{rbp-manni-A-MW-iteration-B}
\mI{Z_r(z)} \babs{B(x)} \dx \, \leq \, c \, r^{\gamma b_{k+1}} \, 
	\big( M^*\big((1+\kabs{Du})^{1+b_k}\big)(z) + F_{k+1}(z) \big)
\end{equation}
for a constant $c$ depending only on $n,L,[u]_{C^{0,\lambda}(Q^+,\R^N)},\alpha,\lambda$ and $\gamma$. We mention that the functions $M^*\big((1+\kabs{Du})^{1+b_k}\big)$ and $F_{k+1}$ belong to the space $L^{s_{k+1}}(Z_{\sigma \rho}(x_0))$, due to the higher integrability of $Du$ and Corollary~\ref{rbp-manni-bem-slice-iteration}, respectively (note $s_{k+1} \in (2, s_k)$).

For the second integral on the right-hand side of \eqref{rbp-manni-A-MW} we argue similarly to above on p. \pageref{rbp-manni-A-inhomogenitaet}: we first assume that we are close to the boundary, i.\,e. $z_n < 2r$. Then, we infer the following estimate from the growth condition (B) on the inhomogeneity, the Caccioppoli inequality (note that $2r \leq \rho \leq \widetilde{\rho}_{cacc}$), the H\"older continuity of $u$ and Poincar\'{e}'s inequality in the boundary version: 
\begin{align}
\label{rbp-manni-A-MW-iteration-b}
r \mI{Z_r(z)} \kabs{b(x,u(x),Du(x))} \dx & \leq \, r \, \mI{Z_r(z)} (L + L_2 \, \kabs{Du}^2) 
\dx \nonumber  \\
	& \leq \, c \, r \, \mI{Z_{2r}(z)} \Big( 1 + \Babs{\frac{u}{r}}^{1+b_k} \, 
	r^{(1-b_k) (\lambda -1)} \Big) \dx \nonumber \\
	& \leq \, c \, r^{1 + (1-b_k)\, (\lambda -1)} 
	\mI{Z_{2r}(z)} \big( 1 + \kabs{Du} \big)^{1 + b_k} \dx \nonumber \\
	& \leq \, c \, r^{b_{k+1}} \, M^*\big((1+\kabs{Du})^{1+b_k}\big)(z) \,,
\end{align}
where in the last line we have employed the fact that $1 + (1-b_k)\, (\lambda -1) \geq  b_{k+1}$
and where the constant $c$ depends only on $n,N,L,L_2,\nu$ and $[u]_{C^{0,\lambda}(Q^+,\R^N)}$. For cylinders in the interior, meaning that $z_n \geq 2r$, we end up with exactly the same estimate using both the Caccioppoli inequality and the Poincar\'{e} inequality with $\kabs{u}$ replaced by $\kabs{u - (u)_{z,2r}}$.

Merging the estimates found in \eqref{rbp-manni-A-MW-iteration-B} and \eqref{rbp-manni-A-MW-iteration-b} together with \eqref{rbp-manni-A-MW} hence yields 
\begin{align*}
\lefteqn{\hspace{-0.75cm} \mI{Z_{\bar{\rho}}(z)} \Babs{\mI{D_{\bar{\rho}}(z')} a_n(y',x_n,u(y',x_n),Du(y',x_n)) \, dy' - 
	\big( a_n(\ccdot,u,Du)\big)_{z,\bar{\rho}} } \dx } 
	 \\
	& = \, \mI{I_{\bar{\rho}}(z_n)} \babs{A_{\bar{\rho}}(x_n) - 
	(A_{\bar{\rho}})_{z_n,\bar{\rho}}} \dx_n \,
	 \leq \, c \, r^{\gamma b_{k+1}} \, 
	\big[ M^*\big((1+\kabs{Du})^{1+b_k}\big)(z) + F_{k+1}(z) \big]
\end{align*}
for a constant $c$ depending only on $n,N,L,L_2,[u]_{C^{0,\lambda}(Q^+,\R^N)},\alpha,\lambda$ and $\gamma$. This is the desired improvement of inequality \eqref{rbp-manni-A-MW-2}. Moreover, $F_{k+1}, M^*\big((1+\kabs{Du})^{1+b_k}\big) \in L^{s_{k+1}}(Z_{\sigma \rho}(x_0))$ holds true. In order to find a fractional Sobolev estimate for the map $x \mapsto a_n(x,u(x),Du(x))$ it still remains to deduce an estimate corresponding to \eqref{rbp-manni-A-MW-3}. To this aim we follow the line of arguments leading to  \eqref{rbp-manni-A-MW-3} and \eqref{u-Du-iteration}: in view of Corollary~\ref{rbp-manni-bem-slice-iteration}, H\"older's inequality and the H\"older continuity of $u$, we see
\begin{align*}
\lefteqn{\hspace{-0.75cm} \mI{Z_{\bar{\rho}}(z)} \Babs{a_n(x,u(x),Du(x)) - \mI{D_{\bar{\rho}}(z')} a_n(y',x_n,u(y',x_n),Du(y',x_n)) \, dy'} \dx } \\
	& \leq \, L \mI{Z_{\bar{\rho}}(z)} \mI{D_{\bar{\rho}}(z')} 
	\babs{Du(x',x_n) - Du(y',x_n)} \dy' \dx \\
	& \qquad {} + 4 \, L \, \mI{Z_{\bar{\rho}}(z)} \mI{D_{\bar{\rho}}(z')} 
	\big( \bar{\rho}^{\, \alpha} + \kabs{u(x',x_n) - u(y',x_n)}^{\alpha} \big) \, 
	\big(1+\kabs{Du(x)}\big) \dy' \dx \\ 
	& \leq \,  c \, r^{\gamma b_{k+1}} \, F_{k+1}(z) 
	+ 4 \, L \, \bar{\rho}^{\, \alpha} \mI{Z_{\bar{\rho}}(z)} \big(1+\kabs{Du(x)}\big) \dx \\
	& \qquad {} + 8 \, L \, \Big( \mI{Z_r(z)}  \kabs{u(x)-(u)_{z,r}}^{\alpha \frac{1+b_k}{b_k}} \dx 
	\Big)^{\frac{b_k}{1 + b_k}}  \Big( \mI{Z_r(z)} \big( 1 + \kabs{Du} )^{1+b_k} \dx 
	\Big)^{\frac{1}{b_k + 1}} \\
	& \leq \, c \, r^{\gamma b_{k+1}} \, 
	\big( M^*\big((1+\kabs{Du})^{1+b_k}\big)(z) + F_{k+1}(z) \big) \,,
\end{align*}
and the constant $c$ depends only on $n,L,[u]_{C^{0,\lambda}(Q^+,\R^N)},\alpha,\lambda$ and $\gamma$. In particular, taking into account  $\bar{\rho} \in [\frac{r}{2},r]$, we infer from the latter two estimates the inequality
\begin{multline*}
\qquad \mI{Z_{r/2}(z)} \babs{a_n(x,u(x),Du(x)) - \big( a_n(\ccdot,u,Du)\big)_{z,r/2} } \dx \\
	\leq \, c \, r^{\gamma b_{k+1}} \, \Big( 
	M^*\big((1+\kabs{Du})^{1+b_k}\big)(z) + F_{k+1}(z) \Big) \,, \qquad
\end{multline*}
where $c$ admits the same dependencies as in the preceding inequalities. Since $F_{k+1}$, $M^*\big((1+\kabs{Du})^{1+b_k}\big)$ belong to $L^{s_{k+1}}(Z_{\sigma \rho}(x_0))$, we may apply the characterization of fractional Sobolev spaces in Lemma~\ref{lemma-frac-char} and Remark~\ref{remark-frac-char}, and we obtain
\begin{equation*}
a_n(\ccdot,u,Du) \in M^{\gamma b_{k+1},s_{k+1}}(Q_{1/(2 \cdot 2^k)}^+,\R^N) \,.
\end{equation*}
Furthermore, there exists $G_{k+1} \in L^{s_{k+1}}(Q_{1/(2^{k+1})}^+,\R^N)$ which satisfies for almost every $x,y \in Q_{1/(2^{k+1})}^+$
\begin{equation*}
\kabs{a_n(x,u(x),Du(x)) - a_n(y,u(y),Du(y))} \, \leq \, 
\kabs{x-y}^{\gamma b_{k+1}} \, \big( G_{k+1}(x) + G_{k+1}(y) \big) \,.
\end{equation*}
We note that $G_{k+1}$ can be calculated from the constant $c$, the functions $M^*((1+\kabs{Du})^{1+b_k})$, $F_{k+1}(z)$ and the restriction on the radius $\rho$ which in turn result in a dependence on the iteration step $k$. For the interior situation we observe that the statements of the Remarks~\ref{bem-manni-innen} remain valid, in particular, the coefficients $a(\cdot,u,Du)$ satisfy a corresponding interior fractional Sobolev estimate.

\subsubsection*{Final conclusion for $\mathbf{Du}$}

Exactly as before on p. \pageref{D_n_u_inherit} the normal derivative $D_n u$ inherits the fractional Sobolev estimate of both the coefficients $a_n(\cdot,u,Du)$ and the tangential derivative $D'u$ (see Corollary~\ref{corollary-tangential-fractional-iteration}). This gives
\begin{equation*}
D u \in M^{\gamma b_{k+1},s_{k+1}}(Q_{1/2^{k+1}}^+,\R^{nN}) \,.
\end{equation*}
At this point we are in the position to use the embeddding
\begin{equation*}
M^{\gamma b_{k+1},s_{k+1}}(Q_{1/2^{k+1}}^+,\R^{nN}) \subset W^{\gamma' \gamma b_{k+1},s_{k+1}}(Q_{1/2^{k+1}}^+,\R^{nN})
\end{equation*}
for all $\gamma' \in (0,1)$. Since $\gamma$ and $\gamma'$ may be chosen arbitrarily close to $1$ (the choice   $\gamma=\gamma'=(\tfrac{n}{n+2\lambda})^{1/2}$ like in the first step is still appropriate for every $k \in \N$), the application of Theorem~\ref{CampInterpolation} yields $Du \in L^{s_{k+1}(1+b_{k+1})}(Q_{1/2^{k+1}}^+,\R^{nN})$. This finishes the iteration and yields:

\begin{lemma}
\label{rbp--manni-fractional-Du}
Let $\alpha \in (0,1)$ and let $u \in W^{1,2}_{\Gamma}(Q_2^+,\R^N) \cap L^{\infty}(Q_2^+,\R^N) \cap C^{0,\lambda}(Q^+,\R^N)$, $\lambda \in (0,1]$, be a weak solution of the Dirichlet problem \eqref{rbp-manni-model} under the assumptions (H1)-(H3) and (B). Then for every $t < \alpha$ there exists $\bar{\rho} = \bar{\rho}(t)$ such that $Du \in W^{t,\bar{s}}(Q_{\bar{\rho}}^+,\R^{nN})$, where $\bar{s}>2$ is a higher integrability exponent depending only $n,N,L,L_2,\nu$ and $[u]_{C^{0,\lambda}(Q^+,\R^N)}$ (but not on $t$).

\begin{proof}
In view of $b_k \nearrow \alpha$, the iteration scheme immediately implies a fractional differentiability result for $Du$: for every $t_1 < \alpha$ there exists $\bar{k} = \bar{k}(t_1)$ such that $Du \in W^{t_1,2}(Q_{1/2^{\bar{k}}}^+,\R^{nN})$. In particular, we may choose $t_1$ sufficiently close to $\alpha$ such that
\begin{equation*}
2 + 2 \alpha \, < \, \frac{2n (1 + t_1)}{n - 2t_1\lambda} \,.
\end{equation*}
In view of Theorem~\ref{CampInterpolation} we hence find $Du \in L^{2 + 2\alpha}$, and at this stage we can indeed stop the iteration: as already observed in Remark~\ref{rem-limit}, in all calculations above the exponent $b_{\bar{k}+1}$ can be replaced by $\alpha$, leading to
\begin{equation*}
D u \in M^{t,s_{\bar{k}+1}}(Q_{1/2^{\bar{k}+1}}^+,\R^{nN})
\end{equation*}
for every $t < \alpha$ and a higher integrability exponent $\bar{s} := s_{\bar{k}+1} > 2$. With the choice $\bar{\rho} := 2^{-\bar{k}-1}$ the proof of the lemma is complete.
\end{proof}
\end{lemma}

\begin{remarkon}
A similar statement was derived for weak solutions to superquadratic nonlinear elliptic systems with inhomogeneities satisfying a controllable growth condition, see \cite[Lemma 6.1]{DUZKRIMIN05}. We easily observe that the method presented in this section does not only apply to inhomogeneities obeying a natural growth condition, but also to those obeying a controllable growth condition. As an advantage of the technique presented here, we note that in the formulation of the previous Lemma~\ref{rbp--manni-fractional-Du} the low dimensional assumption $p > n - 2 - \delta$ for some positive number $\delta$ is not necessary, whereas it was required in the proof of \cite[Lemma 6.1]{DUZKRIMIN05}.
\end{remarkon}

\begin{proof}[of Theorem~\ref{regul-bp-mit-u-nat-2}]
Most of the arguments required here can be recovered from \cite[proof of Theorem 1.1]{DUZKRIMIN05}; for the sake of completeness we sketch briefly the procedure: First, we reduce the general Dirichlet problem \eqref{DP-rbp} with boundary values $u_0$ to the corresponding boundary value problem with zero boundary values, i\,e., $u_0=0$ on $\partial \Omega$ (some attention is needed here: the transformed coefficients still satisfy assumption (H1)-(H3), but the transformed inhomogeneity then satisfies also a slightly different critical growth condition in the sense that 
\begin{equation*}
\kabs{\tilde{b}(x,u,z)} \, \leq \, L + L_2(M) \, \kabs{z}^p
\end{equation*}
for all $(x,u,z) \in \overline{\Omega} \times \R^N \times \R^{nN}$ with $\kabs{u+u_0} \leq M$. Furthermore, the regularity of $\partial \Omega$ allows us to flatten the boundary locally around every boundary point $x_0 \in \partial \Omega$ to end up with a finite number of problems of type \eqref{rbp-manni-model} on cubes. It then suffices to prove that almost every point on $Q^0$ is a regular boundary point, i.\,e. that it belongs to the set $\Reg_{Du}(Q^0)$: since the Hausdorff dimension is invariant under bi-Lipschitz transformations, a standard covering argument then yields the corresponding estimate for the singular boundary points on  $\partial \Omega$, i.\,e. for $\Sing_{Du}(\partial\Omega)$.

In the model situation, \cite[Theorem 1]{ARKHIPOVA03} guarantees that $u$ is H\"older continuous on the regular set $\Reg_u(Q_2^+ \cup \Gamma)$ of $u$ with any exponent $\lambda \in \big(0,1-\frac{n-2}{2}\big)$ and that $\dim_{\Hm}(\Sing_{u}(Q_2^+ \cup Q^0) < n-2$. In particular, the set of singular points is empty if $n=2$. We next observe that the statement in Lemma~\ref{rbp--manni-fractional-Du} still holds true if we replace the cube $Q_1^+$ by any smaller cube $Q_R^+(x_0)$, meaning that we then obtain $Du \in W^{t,\bar{s}}(Q_{\delta R}^+(x_0),\R^{nN})$ for some $\delta(t) > 0$ for all $t < \alpha$ and some $\bar{s}>2$ (independently of the choice of $t$). Therefore, choosing an increasing sequence of sets $B_k \nearrow \Reg_u(Q^+ \cup Q^0)$ with $B_k \subset  \Reg_u(Q^+ \cup Q^0)$ such that $B_k$ is relatively open in $Q^+ \cup Q^0$ for every $k \in \N$, Lemma~\ref{rbp--manni-fractional-Du} allows us to infer that for every $t < \alpha$ and every point $x_0 \in B_k \cap Q^0$ there holds $Du \in W^{t,\bar{s}}(Q_{\delta R}^+(x_0),\R^{nN})$ for some $\delta(t) > 0$. Taking $t \in (2/\bar{s},\alpha)$ and applying the measure density result \cite[Proposition 2.1]{DUZKRIMIN05} (tracing back to Giusti \cite[Proposition 2.7]{GIUSTI03}) we thus find
\begin{equation*}
\dim_{\Hm}\big(\Sing_{Du}(Q^0) \cap Q_{\delta}^+(x_0)\big) \, \leq \, n - \bar{s} t \, < \, n-1 \,.
\end{equation*}
Hence, we find $\dim_{\Hm}(\Sing_{Du}(Q^0) \cap B_k) < n-2t$ for every $k \in \N$ via a covering argument. Keeping in mind $\dim_{\Hm}(\Sing_{u}(Q^+ \cup Q^0) < n-2$, we finally conclude the desired estimate $\dim_{\Hm}(\Sing_{Du}(Q^0)) < n-1$ on the Hausdorff dimension of the singular set for the gradient $Du$ on the boundary. This completes the proof of the main result.
\end{proof}

\begin{remarkon}
We emphasize that the proof also yields a global fractional differentiability result: if $u \in u_0 + W^{1,2}_0(\Omega,\R^N) \cap L^{\infty}(\Omega,\R^N)$ is a weak solutions to \eqref{DP-rbp} under the assumptions of the theorem and with $\alpha \in (0,1)$ arbitrary, then $Du \in W^{t,2}(\Omega,\R^{nN})$ for all $t < \alpha$.
\end{remarkon}


\section{Extensions and open questions related to the dimension reduction}
\label{sec-ext}

In the last section we discuss briefly some extensions and open questions related to regularity results concerning the above regularity theory and to the Hausdorff dimension of the singular set up to the boundary:

\subsection{The nonquadratic case}

In the case $p \in (1,\infty)$ it is convenient to work in terms of the $V$-function (useful algebraic properties are found in \cite[Lemma 2.1]{CAROZZA98} for the subquadratic case, and for the superquadratic case similar inequalities hold true). We start by observing that Mingione \cite{MINGIONE03} was the first to succeed in the dimension reduction for the singular points for $Du$ in the interior of $\Omega$, where $u \in W^{1,p}(\Omega,\R^N) \cap L^{\infty}(\Omega,\R^N)$ is a weak solution (under the standard smallness assumption) to the general system \eqref{DP-rbp} with an inhomogeneity satisfying a natural growth condition and being H\"older continuous with respect to every argument:
\begin{align}
\label{B-cont-1}
\kabs{b(x,u,z) - b(\bar{x},\bar{u},z)} & \leq \, L \, \omega_{\beta_1}\big( \kabs{x - \bar{x}} + \kabs{u - \bar{u}} \big) \, ( 1 + \kabs{z}^p) \\
\label{B-cont-2}
\kabs{b(x,u,z) - b(x,u,\bar{z})} & \leq \, L \, \kabs{z - \bar{z}}^{\beta_2} \, \big(1 + \kabs{z}^2 + \kabs{\bar{z}}^2 \big)^{p-\beta_2} 
\end{align}
for all $x,\bar{x} \in \Omega$, $u, \bar{u} \in \R^N$ and $z, \bar{z} \in \R^{nN}$, $\omega_{\beta_1}(t) = \min\{1,t^{\beta_1}\}$ with $\beta_1,\beta_2 \in (0,1)$ (these condition are in fact only needed to obtain a sharper bound in the superquadratic case). The dimension reduction relies on the fundamental estimate \eqref{rbp-k-start} which iteratively improves the fractional differentiability of $V(Du)$. In a slightly generalized version (including also the subquadratic case) this result can be stated as follows:

\begin{theorem}[cf. \cite{MINGIONE03}, Theorem 2.2]
\label{int-Mingione}
Let $n \leq p + 2$ and $u \in W^{1,p}(\Omega,\R^N) \cap L^{\infty}(\Omega,\R^N)$ be a bounded weak solution to \eqref{DP-rbp} under the conditions (H1)-(H3) and (B), and suppose  $\knorm{u}_{L^{\infty}(\Omega,\R^N)} \leq M$ for some $M>0$ such that $2 \max\{p-1,1\} L_2 M < \nu$. Then we have
\begin{equation*}
\dim_{\Hm}\big(\Sing_{Du}(\Omega)\big) \, \leq \, \max\{ n - \tfrac{p}{p-1}, n - 2 \alpha, n-p \}  \,.
\end{equation*}
If $p>2$ and $\alpha > p / (2p-2)$, then under the assumptions \eqref{B-cont-1}-\eqref{B-cont-2} for $\beta_1 = 2\alpha -1$ and $\beta_2 = (2\alpha - 1)/\alpha$ there holds
\begin{equation*}
\dim_{\Hm}\big(\Sing_{Du}(\Omega)\big) \, \leq \, n - 2 \alpha \,.
\end{equation*}
\end{theorem}

In the general homogeneous case a corresponding interior result \cite{MINGIONE03arch} is valid, and it turned out that the estimates can be extended up to the boundary and that in fact $\Hm^{n-1}$-almost every point in $\overline{\Omega}$ is a regular point for $Du$ independently of the value of $p$ as long as $\alpha > \frac{1}{2}$ and the H\"older continuity of $u$ is known a~priori ($\Hm^{n-1}$-almost everywhere). For this reason one expects that also in the inhomogeneous case under natural growth it should be possible to carry the result of Theorem~\ref{int-Mingione} from the interior up to the boundary. However, apart from some partial result for $p$ close to $2$ (under further restrictions on~$\alpha$), it seems impossible to obtain the boundary regularity statement for all $p \in (1,\infty)$ with the direct approach employed in this paper without any further technical tricks; we emphasize that this problem is not caused by the inhomogeneity but appears also for homogeneous systems where boundary regularity was already proved by the indirect approach. Hence, the direct approach seems to be matched well only in the quadratic situation, and it would be interesting to develop an approach which leads to the known dimension reduction for homogeneous systems, but which can also be applied for inhomogeneous systems under natural growth assumption and a general $p$-growth assumption on the coefficients.

In fact, the only result available so far in this direction is in the two-dimensional case $n=2$, where by means of Morrey-type estimates and a comparison principle the optimal $C^{1,\alpha}_{\rm{loc}}$-regularity of every weak solution is obtained up to the boundary on the regular set of $u$ (the set where $u$ is locally in a neighborhood continuous), and hence outside a set of Hausdorff dimension less than $n-p$:

\begin{theorem}
\label{2d}
Let $p \in (1,\infty)$, $\Omega \subset \R^2$ be a domain of class $C^{1,\alpha}$ and $u_0 \in C^{1,\alpha}(\overline{\Omega},\R^N)$. Let $u \in u_0 + W^{1,p}_0(\Omega,\R^N) \cap L^{\infty}(\Omega,\R^N)$ be a bounded weak solution to \eqref{DP-rbp} under the conditions (H1)-(H3) and (B), and suppose that $\knorm{u}_{L^{\infty}(\Omega,\R^N)} \leq M$ for some $M>0$ such that $2 \max\{p-1,1\} L_2 M < \nu$. Then we have
\begin{equation*}
u \in C^{1,\alpha}(\overline{\Omega},\R^N)
\end{equation*}
in the superquadratic case $p \geq 2$, whereas in the subquadratic case $1 < p < 2$ there hold
\begin{equation*}
u \in C^{1,\alpha}_{\rm{loc}}(\Reg_{u}(\overline{\Omega}),\R^N) \quad \text{ and } \quad \dim_{\Hm}(\overline{\Omega} \setminus \Reg_{u}(\overline{\Omega})) \, < \, n-p \,.
\end{equation*}

\begin{proof}[Sketch]
In the first step of the proof one compares the solution $u$ of the original problem to the solution of the frozen, homogeneous problem, for which good a~priori estimates are known (in particular full $C^1$-regularity). This allows to show that $Du$ belongs to the Morrey-space $L^{p,2-\tau}$ locally on the regular set $\Reg_{u}(\overline{\Omega})$ of $u$ for every $\tau \in (0,2)$, see \cite[p. 2743]{ARKHIPOVA03} if $p \geq 2$ and \cite[p. 317]{BECK09b} if $1<p<2$ for the up-to-the-boundary versions. This decay estimate for the integrals $\int (1 + \kabs{Du})^p \dx$ over balls (possibly intersected with $\Omega$) in terms of the radius is then used to enter again into the comparison argument and to improve the decay of the excess of $Du$ on balls $B_{\rho}(x_0) \subset B_R(x_0)$ contained in the regular set $\Reg_{u}(\overline{\Omega})$ via
\begin{multline*}
\qquad \int_{B_{\rho}(x_0) \cap \Omega} \babs{V(Du)-\big(V(Du)\big)_{B_{\rho}(x_0) \cap \Omega}}^2 \dx \\
	\leq \, 2 \int_{B_{\rho}(x_0) \cap \Omega} \babs{V(Dv)-
	\big(V(Dv)\big)_{B_{R}(x_0) \cap \Omega}}^2 \dx + 
	2 \int_{B_{\rho}(x_0) \cap \Omega} \babs{V(Dv)-V(Du)}^2 \dx \,, \qquad
\end{multline*}
where $v \in u + W^{1,p}(B_{R}(x_0) \cap \Omega,\R^N)$ is the weak solution to $\diverg a(x_0,(u)_{B_{R}(x_0) \cap \Omega},Dv) = 0$. Taking into account the Morrey regularity and proceeding as in \cite[Chapter 9]{KRIMIN05}, we then estimate the first integral on the right-hand side via the decay properties for the comparison map (see \cite[3.12]{CAMPANATO87} and \cite[(4.34)]{BECK08}) combined with a good choice of the radius $R$ as a power of $\rho$. The second integral is under control via an estimate similar to the ones of \cite[(34)]{ARKHIPOVA03} or \cite[(4.16)]{BECK09b}, with the difference that the terms involving $s$ and $\delta$, respectively, do not appear if also the higher integrability of the comparison solution is kept in mind. As a consequence, the reasoning in \cite{KRIMIN05} applies, and a sufficiently small choice of $\tau$ shows that the previous excess integral is bounded by $\rho^{2+\sigma}$ for some $\sigma>0$. This corresponds to H\"older continuity of $V(Du)$ and therefore of $Du$ with some small exponent, which is then improved to the optimal one by standard regularity theory.
\end{proof}
\end{theorem}

\subsection{Systems with coefficients $\mathbf{a(x,z)}$}

In the introduction we already spent some words on the situation where the coefficients do not explicitly depend on the weak solution $u$, but only on its gradient and the independent variable. In the interior it turned out \cite{MINGIONE03arch,MINGIONE03} that the assumption $n \leq p+2$ (or priori H\"older continuity of $u$) is no longer needed (because -- roughly speaking~-- possible singularities of $u$ do propagate to the coefficients only in the gradient variable). For this reason the dimension reduction follows in only one step, and it is further valid without any restrictions on the space dimension $n$. In the indirect approach of \cite{DUZKRIMIN05} the same reasoning was applied in order to extend these interior results up to the boundary. As a consequence, in case of homogeneous systems there holds:

\begin{theorem}[\cite{DUZKRIMIN05}]
\label{regul-bp-ohne-u}
Consider $p \in (1,\infty)$, $\alpha > \frac{1}{2}$. Let $\Omega$ be a domain of class $C^{1,\alpha}$ and $u_0 \in C^{1,\alpha}(\overline{\Omega},\R^N)$. Assume that $u \in u_0 + W^{1,p}_0(\Omega,\R^N)$ is a weak solution to
$\diverg a(x,Du) = 0$ in $\Omega$ under the assumptions (H1)-(H3). Then $\Hm^{n-1}$-almost every boundary point is a regular point for $Du$.
\end{theorem}

It is not clear whether the result of Theorem~\ref{regul-bp-mit-u-nat-2} can be improved to such vector fields which do not explicitly depend on $u$, in the sense that the existence of regular boundary points in that case is valid for all dimensions $n \geq 2$. The first problem arises in the preliminary estimate \eqref{rbp-k-start}, where in general no positive power of $h$ for the last integral -- which came from $|\int b(x,u,Du) \tau_{e,-h}(\eta^2 \tau_{e,h} u) dx|$~-- can be produced. Hence, to obtain some fractional differentiability of the system, additional regularity assumptions on the inhomogeneity are required (such as conditions \eqref{B-cont-1}-\eqref{B-cont-2} of H\"older continuity with respect to all variables), which then allow to use the formula for partial integration for finite differences also in the integral involving the right-hand side. In the interior this gives $\dim_{\Hm}(\Sing_{Du}(\Omega)) \leq n-2\alpha$ in all dimensions. However, in the direct approach presented in this paper, these assumptions do not seem to lead to fractional Sobolev estimates for $x \mapsto a_n(x,Du)$, see e.g. the derivation of \eqref{rbp-manni-A-inhomogenitaet} where the H\"older continuity was the crucial ingredient.

\subsection{Optimality of the Hausdorff dimension}

It is not clear to what extent the estimates for the Hausdorff dimension of the singular set may still be improved, neither in our main statements nor in the cited results \cite{CAMPANATO87,ARKHIPOVA03,BECK09b,MINGIONE03arch,MINGIONE03,DUZKRIMIN05} on the Hausdorff dimension of the singular set of $u$ and $Du$. Up to now, the bound for the singular set of $Du$ depends on the parameter $\alpha$. While one cannot rule out that the dependence on $\alpha$ is only due to technique, it is believed that this dependence is a structural feature of the problem concerning the Hausdorff dimension of the singular set. However, the literature lacks appropriate counterexamples. As a consequence, the question of the existence of regular boundary points for H\"older exponents $\alpha \in (0,\frac{1}{2})$ remains open for general nonlinear systems of the type considered above. Recently, it was observed by Kristensen and Mingione~\cite{KRIMIN08} that the H\"older continuity assumption in $x$ can be relaxed to a fractional Sobolev dependence. More precisely, still assuming the H\"older continuity assumption (H3) with now an arbitrarily small $\alpha>0$ (as a consequence, we have $C^{1,\alpha}_{\rm{loc}}$ regularity on the regular set and the characterization of the singular set remains unchanged), we further assume for some $\beta>0$:
\begin{align*}
\text{{\bf (H4)}} \qquad & \text{There exists a function } g \in L^{\upsilon}(\Omega) \text{ with } g \geq 0 \text{ and } \upsilon > (1+\beta)/\beta \text{ such that }\\
	& \hspace{1cm} \kabs{a(x,u,z) -  a(\bar{x}, \bar{u}, z)} \, \leq \, L 
	\, \big(1 + \kabs{z}^2\big)^{\frac{p-1}{2}} \, \big[ \kabs{x-\bar{x}}^{\beta} (g(x) + g\big(\bar{x})\big) 
	+ \omega_{\beta}\big(\kabs{u-\bar{u}}\big) \big]
\end{align*}
for almost all $x,\bar{x} \in \Omega$ and all $u,\bar{u} \in \R^N$, $z \in R^{nN}$ (and $\omega_{\beta}(s) \leq \min\{1,s^{\alpha}\}$). The function $g$ again plays the role of a fractional derivative, see Section~\ref{fractSob}. Obviously, condition (H4) is weaker than (H3) in the case $\alpha = \beta$ and $\upsilon < \infty$, and it actually turns out that -- independently of the value of~$\alpha$~-- the result of Theorem~\ref{regul-bp-mit-u-nat-2} still holds true, provided that $\beta > \frac{1}{2}$. In particular, we still get existence of regular boundary points (even though in a regular point the exponent of H\"older continuity of $Du$ is only $\alpha$ and in general not $\beta$):

\begin{theorem}
\label{thm-frac}
Consider $n \in \{2,3,4\}$, $\alpha,\beta \in (0,1]$ with $\beta \geq 1/2$. Let $\Omega \subset \R^n$ be a domain of class $C^{1,\alpha}$ and $u_0 \in C^{1,\alpha}(\overline{\Omega},\R^N)$. Assume further that $u \in u_0 + W^{1,2}_0(\Omega,\R^N) \cap L^{\infty}(\Omega,\R^N)$ is a weak solution of the Dirichlet problem \eqref{DP-rbp} under the assumptions (H1)-(H4) and (B), and suppose that $\knorm{u}_{L^{\infty}(\Omega,\R^N)} \leq M$ for some $M>0$ such that $2 L_2 M < \nu$. Then $\Hm^{n-1}$-almost every boundary point is a regular point for $Du$.

\begin{proof}
The strategy of the proof of the theorem is the same as the one of Theorem~\ref{regul-bp-mit-u-nat-2}, and we immediately get into the study of the transformed system \eqref{rbp-manni-model} under the assumption of a~priori H\"older continuity of $u$ with exponent $\lambda$.  We first infer $Du \in L^{2 + 2 \alpha}(Q^+_{\bar{\rho}},\R^{nN})$ from Lemma~\ref{rbp--manni-fractional-Du} from a first iteration using only assumptions (H1)-(H3) and (B) (alternatively we can use a simple higher integrability result via Gehring's Lemma), and now start a new iteration as in Section~\ref{sec-iteration} by taking into consideration the additional assumption (H4): for this purpose we define a sequence $(\beta_k)_{k \in \N}$
\begin{equation*}
\beta_0 \, := \, \frac{\alpha}{2}, \qquad \beta_{k+1} \, := \, \beta_k + 
\min \Big \{ \frac{\lambda}{2}, \frac{\beta_k}{2(1+\beta_k)} \Big \} \, ( \beta - \beta_k ) \,,
\end{equation*}
and we observe that it is bounded and strictly increasing with limit $\beta$. Observing that only the fractional dependence of the coefficients in the $x$-variable has changed, we now have to re-estimate the terms involving differences of the coefficients with respect to the $x$-variable, and then all the statement of Section~\ref{sec-iteration} remain true for $b_k$ replaced by $\beta_k$ (and on smaller half-cubes). In fact, there are only two new terms. Under the H\"older continuity assumption (H3), these were estimated trivially, but they now need to be investigated more carefully: the first occurs in the proof of Proposition~\ref{rbp-manni-diff-higher-iteration} (assuming that $Du$ is integrable to a power greater than $2 + 2\beta_k$), when we derive a suitable substitute for the preliminary estimate \eqref{rbp-k-start} to find the inequality corresponding to \eqref{rbp-manni-diff-iteration}. Actually, only the integral involving ${\cal A}(h)$ needs to be adjusted: Keeping in mind the integrability assumption $g \in L^{\upsilon}(Q_2^+) \subset L^{(1+\beta)/\beta}(Q_2^+)$ and the H\"older continuity of $u$ we calculate with (H3), (H4), Young's and H\"older's inequality:
\begin{align*}
\lefteqn{ \hspace{-0.5cm} \int_{Q^+} \kabs{{\cal A}(h)} \, \kabs{D(\eta^2 \tau_{e,h} u )} \dx
	\, \leq \, L \, \kabs{h}^{\beta \frac{1+\beta}{\beta} \frac{\beta_k}{1+\beta_k}} } \\
	& \qquad {} \times \int_{Q^+} \big(1 + \kabs{Du(x+he)}\big) \, 
	\big( g(x+he) + g(x)\big)^{\frac{1+\beta}{\beta} \frac{\beta_k}{1+\beta_k}} 
	\, \big( \eta^2 \kabs{\tau_{e,h} Du} 
	+ 2 \eta \kabs{D\eta} \kabs{\tau_{e,h}u }\big) \dx \\
	& \leq \, \epsilon \int_{Q^+} \eta^2 \kabs{\tau_{e,h} Du}^2 \dx 
	 + \Big[ c(\epsilon,L) \, 
	\kabs{h}^{2\beta \frac{1+\beta}{\beta} \frac{\beta_k}{1+\beta_k}} 
	+ c(L,\knorm{D\eta}_{L^\infty}) \, 
	\kabs{h}^{\beta \frac{1+\beta}{\beta} \frac{\beta_k}{1+\beta_k} + 1} \Big] \\
	& \qquad {} \times
	\Big( \int_{Q^+ \cap \supp(\eta)} \big(1 + \kabs{Du(x+he)}\big)^{2 + 2 \beta_k} \dx
	\Big)^{\frac{1}{1+\beta_k}}
	\Big( \int_{Q^+} \kabs{g}^{\frac{1+\beta}{\beta}} \dx \Big)^{\frac{\beta_k}{1+\beta_k}} \,,
\end{align*}
where we also have used standard estimates for finite differences (note that the exponent $\alpha$ was for simplicity treated as 0 in the powers of $\kabs{h}$). To conclude \eqref{rbp-manni-diff-iteration} it then suffices to observe that both powers of $\kabs{h}$ are at least $2\beta_{k+1}$ and that the other terms are estimates exactly as before (but using (H4) instead of (H3) for estimating differences of the coefficients with respect to the $u$-variable).

The second new term arises in the fractional Sobolev estimate for $a_n(x,u,Du)$ (and in turn in the same way also for $D_nu$): it occurs for the first time in the estimate for $B(x)$ and can be dealt with as follows:
\begin{multline*}
\quad \mI{Z_r(z)} \babs{a(x,u(x),Du(x)) - a(z,(u)_{z,r},Du(x))} \dx \\ 
	\leq c(n) \, r^{\beta \frac{1+\beta}{2\beta} \frac{\beta_k}{1+\beta_k}} 
	\Big[ \mI{Z_r(z)} \big(1 + \kabs{Du(x+he)}\big)^{1 + \beta_k} \dx + 
	 \mI{Z_r(z)} \big( g(z) + g(x) \big)^{\frac{1 + \beta}{2 \beta}} \dx \Big] \,. \quad
\end{multline*}
The right-hand side is then estimated further via the maximal function.

With these two adjustments, the proof of the theorem then continues as before, leading to the existence of regular boundary points for $\alpha \in (0,1)$ arbitrarily, provided that $\beta \geq \frac{1}{2}$.
\end{proof}
\end{theorem}


\footnotesize
\begin{thebibliography}{10}

\bibitem{ACEFUS89}
E.~Acerbi and N.~Fusco, \emph{{Regularity for minimizers of non-quadratic
  functionals: the case $1<p<2$}}, J. Math. Anal. Appl. \textbf{{\bf 140}}
  (1989), 115--135.

\bibitem{ADAMS75}
R.~A. Adams, \emph{{Sobolev Spaces}}, Academic Press, New York, 1975.

\bibitem{ARKHIPOVA96}
A.~Arkhipova, \emph{On the regularity of the solutions of boundary-value
  problem for quasilinear elliptic systems with quadratic nonlinearity}, J.
  Math. Sci., New York \textbf{{\bf 80}} (1996), no.~6, 2208--2225.

\bibitem{ARKHIPOVA97}
A.~Arkhipova, \emph{{On the Neumann problem for nonlinear elliptic systems with
  quadratic nonlinearity}}, St. Petersbg. Math. J. \textbf{{\bf 8}} (1997),
  no.~5, 845--877.

\bibitem{ARKHIPOVA03}
A.~Arkhipova, \emph{{Partial regularity up to the boundary of weak solutions of
  elliptic systems with nonlinearity $q$ greater than two}}, J. Math. Sci.
  (N.\,Y.) \textbf{{\bf 115}} (2003), 2735--2746.

\bibitem{BECK07}
L.~Beck, \emph{Partial regularity for weak solutions of nonlinear elliptic
  systems: the subquadratic case}, Manuscr. Math. \textbf{{\bf123}} (2007),
  no.~4, 453--491.

\bibitem{BECK08}
L.~Beck, \emph{Boundary regularity results for weak solutions of subquadratic
  elliptic systems}, Ph.D. thesis, Universit\"at Erlangen-N\"urnberg, 2008.

\bibitem{BECK09b}
L.~Beck, \emph{{Partial H\"older continuity for solutions of subquadratic
  elliptic systems in low dimensions}}, J. Math. Anal. Appl. \textbf{{\bf354}}
  (2009), no.~1, 301--318.

\bibitem{CAMPANATO82a}
S.~Campanato, \emph{{Differentiability of the solutions of nonlinear elliptic
  systems with natural growth}}, Ann. Mat. Pura Appl. Ser. 4 \textbf{{\bf 131}}
  (1982), 75--106.

\bibitem{CAMPANATO82}
S.~Campanato, \emph{{H\"older continuity and partial H\"older continuity results for
  $W^{1,q}$-solutions of non-linear elliptic systems with controlled growth}},
  Rend. Sem. Mat. Fis. Milano \textbf{{\bf 52}} (1982), 435--472.

\bibitem{CAMPANATO87}
S.~Campanato, \emph{{Elliptic systems with non-linearity $q$ greater or equal to
  two. Regularity of the solution of the Dirichlet problem}}, Ann. Mat. Pura
  Appl. Ser. 4 \textbf{{\bf 147}} (1987), 117--150.

\bibitem{CAMPCAN81}
S.~Campanato and P.~Cannarsa, \emph{{Differentiability and partial H\"older
  continuity of the solutions of non-linear elliptic systems of order $2m$ with
  quadratic growth.}}, Ann. Sc. Norm. Super. Pisa Ser. IV \textbf{{\bf 8}}
  (1981), 285--309.

\bibitem{CAROZZA98}
M.~Carozza, N.~Fusco, and G.~Mingione, \emph{{Partial Regularity of Minimizers
  of Quasiconvex Integrals with Subquadratic Growth}}, Ann. Mat. Pura Appl.
  Ser. 4 \textbf{{\bf 175}} (1998), 141--164.

\bibitem{COLOMBINI71}
F.~Colombini, \emph{{Un teorema di regolarit\`{a} alla frontiera per soluzioni
  di sistemi ellittici quasi lineari}}, Ann. Sc. Norm. Super. Pisa Ser. III
  \textbf{{\bf 25}} (1971), 15--161.

\bibitem{DEGIORGI68}
E.~De~Giorgi, \emph{{Un esempio di estremali discontinue per un problema
  variazionale di tipo ellittico}}, Boll. Unione Mat. Ital., IV. \textbf{{\bf
  1}} (1968), 135--137.

\bibitem{DOMOKOS04}
A.~Domokos, \emph{{On the regularity of $p$-harmonic functions in the
  Heisenberg group}}, Ph.D. thesis, University of Pittsburgh, 2004.

\bibitem{DUGROKRO04}
F.~Duzaar, J.~F. Grotowski, and M.~Kronz, \emph{{Partial and full boundary
  regularity for minimizers of functionals with nonquadratic growth}}, J.
  Convex Anal. \textbf{{\bf 11}} (2004), no.~2, 437--476.

\bibitem{DUZKRIMIN05}
F.~Duzaar, J.~Kristensen, and G.~Mingione, \emph{{The existence of regular
  boundary points for non-linear elliptic systems}}, {J. Reine Angew. Math.}
  \textbf{{\bf 602}} (2007), 17--58.

\bibitem{DUMIN05}
F.~Duzaar and G.~Mingione, \emph{Second order parabolic systems, optimal
  regularity, and singular sets of solutions}, Ann. Inst. Henri Poincar\'e
  Anal. Non Lin\'eaire \textbf{{\bf 22}} (2005), no.~6, 705--751.

\bibitem{GIAQUINTA78}
M.~Giaquinta, \emph{{A counterexample to the boundary regularity of solutions
  to elliptic quasilinear systems}}, Manuscr. Math. \textbf{{\bf 24}} (1978),
  217--220.

\bibitem{GIAMOD79}
M.~Giaquinta and G.~Modica, \emph{{Almost-everywhere regularity results for
  solutions of non linear elliptic systems}}, Manuscr. Math. \textbf{{\bf 28}}
  (1979), 109--158.

\bibitem{GIUSTI03}
E.~Giusti, \emph{{Direct Methods in the Calculus of Variation}}, World
  Scientific Publishing, Singapore, 2003.

\bibitem{GIUSMIRA68}
E.~Giusti and M.~Miranda, \emph{{Sulla Regolarit\`{a} delle Soluzioni Deboli di
  una Classe di Sistemi Ellitici Quasi-lineari}}, Arch. Rational Mech. Anal.
  \textbf{{\bf 31}} (1968), 173--184.

\bibitem{GIUSMIRA68EX}
E.~Giusti and M.~Miranda, \emph{{Un esempio di soluzioni discontinue per un problema di minimo
  relativo ad un integrale regolare del calcolo delle variazioni}}, Boll.
  Unione Mat. Ital., IV. Ser. \textbf{{\bf 1}} (1968), 219--226.

\bibitem{GROTOWSKI02QL}
J.~F. Grotowski, \emph{{Boundary regularity for quasilinear elliptic systems}},
  Commun. Partial Differ. Equations \textbf{{\bf 27}} (2002), no.~11-12,
  2491--2512.

\bibitem{GROTOWSKI02}
J.~F. Grotowski, \emph{{Boundary regularity results for nonlinear elliptic systems}},
  Calc. Var. Partial Differ. Equ. \textbf{{\bf 15}} (2002), 353--388.

\bibitem{HAJLASZ95}
P.~Haj{\l}asz, \emph{{Geometric approach to Sobolev spaces and badly
  degenerated elliptic equations}}, Proceedings of the Banach Center
  minisemester on nonlinear analysis and applications, \emph{GAKUTO Int. Ser.,
  Math. Sci. Appl.} \textbf{{\bf7}} (1995), 141--168.

\bibitem{HAJLASZ96}
P.~Haj{\l}asz, \emph{{Sobolev spaces on an arbitrary metric space}}, Potential Anal.
  \textbf{{\bf5}} (1996), no.~4, 403--415.

\bibitem{HAJKOS00}
P.~Haj{\l}asz and P.~Koskela, \emph{{Sobolev met Poincar\'{e}}}, Mem. Am. Math.
  Soc. \textbf{{\bf688}} (2000), 101 p.

\bibitem{HAMBURGER07}
C.~Hamburger, \emph{{Partial boundary regularity of solutions of nonlinear
  superelliptic systems}}, Boll. Unione Mat. Ital., Sez. B, Artic. Ric. Mat.
  (8) \textbf{{\bf10}} (2007), no.~1, 63--81.

\bibitem{IVERT79}
P.-A. Ivert, \emph{{Regularit\"atsuntersuchungen von L\"osungen elliptischer
  Systeme von quasilinearen Differen\-tialgleichungen zweiter Ordnung}},
  Manuscr. Math. \textbf{{\bf 30}} (1979), 53--88.

\bibitem{KRIMIN05a}
J.~Kristensen and G.~Mingione, \emph{{The Singular Set of $\omega$-minima}},
  Arch. Rational Mech. Anal. \textbf{{\bf 177}} (2005), 93--114.

\bibitem{KRIMIN05}
J.~Kristensen and G.~Mingione, \emph{{The Singular Set of Minima of Integral Functionals}}, Arch.
  Rational Mech. Anal. \textbf{{\bf 180}} (2006), no.~3, 331--398.

\bibitem{KRIMIN08}
J.~Kristensen and G.~Mingione, \emph{Boundary regularity in variational problems}, Arch. Rational
  Mech. Anal. (to appear).

\bibitem{KRONZHABIL}
M.~Kronz, Habilitationsschrift, Erlangen, in preparation.

\bibitem{MINGIONE03}
G.~Mingione, \emph{Bounds for the singular set of solutions to non linear
  elliptic systems}, Calc. Var. Partial Differ. Equ. \textbf{{\bf 18}} (2003),
  no.~4, 373--400.

\bibitem{MINGIONE03arch}
G.~Mingione, \emph{{The Singular Set of Solutions to Non-Differentiable Elliptic
  Systems}}, Arch. Rational Mech. Anal. \textbf{{\bf 166}} (2003), 287--301.

\bibitem{PEPE71}
L.~Pepe, \emph{{Risultati di regolarit\`{a} parziale per le soluzioni
  $\mathcal{H}^{1,p}$ $1 < p < 2$ di sistemi ellittici quasi lineari}}, Ann.
  Univ. Ferrara, N. Ser., Sez. VII \textbf{{\bf 16}} (1971), 129--148.

\bibitem{SIMON96}
L.~Simon, \emph{{Theorems on Regularity and Singularity of Energy Minimizing
  Maps}}, Birk\-h\"auser-Verlag, Basel-Boston-Berlin, 1996.

\bibitem{TOLK83}
P.~Tolksdorf, \emph{Everywhere-regularity for some quasilinear systems with a
  lack of ellipticity}, Ann. Mat. Pura Appl., IV. Ser. \textbf{{\bf 134}}
  (1983), 241--266.

\bibitem{UHLENBECK77}
K.~Uhlenbeck, \emph{{Regularity for a class of nonlinear elliptic systems}},
  Acta Math. \textbf{{\bf138}} (1977), 219--240.

\bibitem{WIEGNER76}
M.~Wiegner, \emph{{A-priori Schranken f\"ur L\"osungen gewisser elliptischer
  Systeme}}, Manuscr. Math. \textbf{{\bf 18}} (1976), 279--297.

\end{thebibliography}
\end{document}